\documentclass[mathpazo]{cicp}
\usepackage{algorithmic}
\usepackage{booktabs}
\usepackage{listings}
\usepackage{mathtools} 
\usepackage{mathalfa}
\usepackage{tikz}
\usepackage{pgfplots}
\pgfplotsset{compat=1.17} 

\usepackage{soul} 
\setstcolor{red}

\begin{document}
\title{A Stochastic Tamed-Regularized Random Batch Method for Uncertain Particle Systems with Singularities}

\author[W.~Wang, L.~Li, and Y.~Lin]{Wenyi Wang\affil{1}\comma, Lei Li\affil{1,2}, and Yiwen Lin\affil{1}\corrauth}
\address{\affilnum{1}\ School of Mathematical Sciences, Shanghai Jiao Tong University, Shanghai 200240, China. \\
         \affilnum{2}\ Institute of Natural Sciences, MOE-LSC, Shanghai Jiao Tong University, Shanghai 200240, China}
\emails{{\tt wangwenyi2004@sjtu.edu.cn} (W.~Wang), {\tt leili2010@sjtu.edu.cn} (L.~Li), {\tt linyiwen@sjtu.edu.cn} (Y.~Lin)}


\begin{abstract}
The coefficients of particle systems often exhibit uncertainties due to incomplete knowledge. For particle systems with uncertainties and singularities in the drift and diffusion terms, existing methods may encounter numerical blow-up. In this paper, we propose a novel stochastic tamed-regularized Random Batch (str-RB) method for such uncertain particle systems with singularities. Specifically, since singularities in the diffusion term are more challenging than those in the drift term, building on a tamed scheme, we design an exponential regularization strategy that more effectively mitigates singularities in the diffusion term. This strategy ensures that the stochastic Galerkin (sG) projection coefficients remain finite even when the corresponding integrals for the original singular diffusion coefficient diverge. The sG method is applied to handle the random coefficients, and a random batch approach is then employed to reduce the computational complexity of large-scale simulations. Several theoretical properties of the proposed algorithm are established. Numerical tests demonstrate that the str-RB algorithm exhibits weak convergence, achieves lower errors than major existing methods across singular tests, and remains stable in a wider range of singular regimes. This work provides a robust and accurate approach for uncertain particle systems with singularities, laying a foundation for high-performance solvers in uncertain many-body dynamics.
\end{abstract}


\maketitle
\section{Introduction}
\par Particle systems are fundamental mathematical models for describing the collective dynamics of large numbers of interacting agents, and they play an indispensable role in modern scientific computing. With wide applications  in computational physics \cite{27}, chemistry \cite{28}, materials science \cite{29}, and biological engineering \cite{30}, these models effectively capture complex many-body phenomena. Particle systems with singular coefficients, in particular, arise in diverse applications such as flocking models \cite{32} and semiconductor device simulation \cite{33}. In this study, we design a new computational method for a class of uncertain particle systems that exhibit singularities. 
\par Most existing studies on particle systems concentrate on models with deterministic drift and diffusion coefficients. In practical scientific and engineering applications, however, the coefficients often involve uncertainties arising from incomplete knowledge of the system parameters. This necessitates the development of numerical methods for particle systems with random coefficients. The Stochastic Galerkin (sG) method is a fundamental and widely adopted  approach for uncertainty quantification, which transforms differential equations with random coefficients into systems with deterministic coefficients. In \cite{14,15,17,19}, the sG method was implemented for diffusion equations, the Boltzmann equation, Hamilton--Jacobi equations, and Navier--Stokes equations. The method has also been applied to other elliptic equations, hyperbolic equations, and particle systems \cite{16,18,22}. Nevertheless, a direct application of the sG method to large-scale particle systems incurs a high computational cost, which motivates the use of a random batch strategy to reduce complexity.
\par The Random Batch Method (RBM), first proposed in \cite{1}, is an efficient approach for particle systems that significantly reduces the computational cost. Its core idea is to randomly divide the particles into multiple small batches at each time step, so that each particle interacts only with particles within the same batch. From the theoretical perspective, rigorous proofs for error estimates and convergence properties for the RBM and some of its variants were established in \cite{1,12,13}. 
The RBM has been applied to a variety of classical and quantum systems, including the Poisson-Nernst-Planck system, the Poisson-Boltzmann system, and systems driven by L\'{e}vy noise \cite{4,5,6,7}. Moreover, it has been integrated with the Ewald summation \cite{8,9} for Coulomb interactions, splitting strategies for singular kernels \cite{10}, and the sum-of-Gaussians decomposition for long-range interactions \cite{11} to develop composite methods. For particle systems with uncertainties, a particle sG method was proposed in \cite{3} by combining a random batch approach with the sG method. This method is effective for uncertain particle systems with smooth coefficients and serves as one of the baselines in the numerical tests of this work. 
\par While the existing particle sG method performs well when the drift and diffusion terms depend smoothly on random variables, it tends to become numerically unstable when the coefficients have singularities of higher order. In particular, in combination with Brownian motion, such singularities in the diffusion term may introduce large errors. Although tamed schemes \cite{20,21,23,24} and some extensions \cite{2,25,26} have been developed, these algorithms are designed primarily for the drift term and may produce relatively large errors induced by singularities in the diffusion term. To further mitigate these errors, we propose a novel stochastic tamed-regularized Random Batch (str-RB) algorithm, which simultaneously tames the singular drift term and regularizes the singular diffusion term. Numerical experiments demonstrate weak convergence and its ability to yield more accurate simulations for uncertain particle systems with singularities. 
\par The proposed str-RB algorithm not only tames the drift term but, more critically, introduces a tailored exponential regularization strategy for the singular diffusion term. Firstly, we construct a regularization function for the diffusion term, which ensures that the sG projection can be performed stably, and the random-coefficient system is then converted into a stochastic differential equation (SDE) with deterministic coefficients. The dimension of this SDE is proportional to the number of particles, so solving it directly entails a computational cost that scales quadratically with the particle number. To alleviate this high computational cost, a variant of the RBM is applied to the SDE, which reduces the computational complexity to the same order as the particle number. Subsequently, to avoid numerical instability caused by singularities in the drift term, the modified tamed scheme in \cite{2} is employed to solve the SDE with deterministic coefficients. It is noted that the exponential regularization strategy plays a pivotal role in  mitigating the instability caused by highly singular diffusion terms, and several theoretical guarantees of this strategy, including numerical stability and an analysis of the regularization-induced model perturbation, are established to demonstrate its effect. Finally, we conduct a range of numerical tests, focusing primarily on cases where the diffusion term exhibits challenging singularities, to validate the convergence and effectiveness of the str-RB algorithm. The numerical results demonstrate weak convergence and show that, in challenging singular cases, the algorithm yields lower errors than both existing methods and the enhanced baseline, with error reductions of several orders of magnitude in some cases, while achieving comparable or better performance on smooth and low‑order singular problems. This demonstrates that the str-RB algorithm handles a broader class of singular problems that are difficult for existing methods, and remains competitive in easier cases, maintaining both numerical stability and high accuracy.
\par The paper is organized as follows. In Section 2 we describe the particle systems considered in this work. The key ingredients and the complete procedure of the str-RB algorithm are presented in Section 3. Section 4 discusses the theoretical properties of the str-RB algorithm. Section 5 presents numerical tests that validate the convergence, numerical stability and effectiveness of the str-RB algorithm. The paper is concluded in Section 6.
\section{Particle systems}
\par We first introduce the random-coefficient particle systems considered in this study, whose mean-field limit is the Vlasov-Fokker-Planck equation with random coefficients. These random particle systems take the following form:
\begin{equation}
\left\{
\begin{aligned}
dX_i(z, t) &= V_i(z, t)\, dt,\\
dV_i(z, t) &= \frac{1}{N} \sum_{j=1}^N P(z, X_i, X_j)\,(V_j - V_i)\, dt + \sqrt{2D(z)}\, dW_i,
\end{aligned}
\right.
\qquad i=1,\dots, N,
\label{2}
\end{equation}
where $N$ is the number of particles, $z$ is a random variable modeling the uncertainties in the coefficients, $V_i(z,t)$ and $X_i(z,t)$ represent the velocity and position of the $i$-th particle at time $t$, $P(z, X_i, X_j)$ is an interaction function characterizing the interaction between the $i$-th and $j$-th particles, $D(z)$ is a positive diffusion coefficient, and  $\{W_i\}_{i=1}^N$ are independent Brownian motions.
\par By taking $N \to \infty$, one derives the mean-field Vlasov-Fokker-Planck equation with uncertainties:
\begin{align*}
\partial_t f(x,v,t,z) + v \cdot \nabla_x f(x,v,t,z) =  Q(f,f)(x,v,t,z),
\end{align*}
where
\begin{align*}
 Q(f,f) = \nabla_v \cdot \left( \mathcal{P}[f] f + \nabla_v (D(z) f) \right) ,  
\end{align*}
and
\begin{align*}
    \mathcal{P}[f](x, v, t, z) = \int\nolimits_{\mathbb{R}^2} P(z, x, x_*) (v - v_*) f(x_*, v_*, t, z) \, dv_* dx_*.
\end{align*}
\par Efficient numerical methods have been well developed for particle systems with smooth drift and diffusion terms \cite{3}. However, singularities with respect to the random variable in both drift and diffusion terms bring new computational difficulties. This work is devoted to developing numerical algorithms for random-coefficient particle systems of the form Eq.~\eqref{2} with singular coefficients, as detailed in Section 3.
\section{The str-RB algorithm}
\par To address the numerical instability caused by singularities in the diffusion term of Eq. \eqref{2}, we propose an exponential regularization strategy. A random batch approach and a modified tamed scheme are then adopted to reduce the computational complexity and to tame the singularities in the drift term, respectively. The resulting algorithm is referred to as the \textit{stochastic tamed-regularized Random Batch} (str-RB) algorithm. This section details the key components of the str-RB algorithm and presents the complete procedure. Specifically, in Section 3.1, we introduce the exponential regularization strategy designed for singular diffusion terms.  Section 3.2 presents the sG approximation of the regularized particle system. The resulting deterministic-coefficient SDE is then solved by a random batch approach together with the modified tamed scheme \cite{2}, which will be introduced in Sections 3.3 and 3.4, respectively. The complete algorithm is summarized in Section 3.5.
\subsection{Exponential regularization strategy}
\par The singularities in the diffusion term $D(z)$, in combination with Brownian motion, tend to cause more severe numerical instability than those in the drift term. It is therefore necessary to regularize $D(z)$ before performing the sG projection, so that the sG projection can be carried out stably. To this end, we propose an exponential regularization strategy, in which a regularization function $\psi(r)$ is designed to provide stronger suppression of large values of $D(z)$ than algebraic damping \cite{2,21}. To be specific, the exponential regularization function is defined as
\begin{align}
\psi(r) = \exp\left( -\tau^{\alpha_{\text{diff}}}  \gamma_{\text{diff}}  r^{\beta_{\text{diff}}} \right)\label{3-4}
\end{align}
where $\tau$ is the time step, $\beta_{\text{diff}}\in(0,1)$, and $\alpha_{\text{diff}} ~\text{and}~\gamma_{\text{diff}}$ are non-negative parameters, each playing a distinct role in the regularization. Specifically, 
$\gamma_{\text{diff}}$ sets the fixed regularization strength; $\alpha_{\text{diff}}$ controls the rate at which the regularization weakens as the time step $\tau$ decreases; $\beta_{\text{diff}}$ determines how rapidly the regularization strengthens as the diffusion coefficient $D(z)$ grows.  In this exponential regularization strategy, the regularized diffusion term
\begin{align}
\displaystyle\tilde{D}(z)=\Psi (D(z))=D(z)\cdot \psi(D(z)) \label{3-3}   
\end{align}
is adopted in place of the original $D(z)$, thereby mitigating the singularity-induced instability and effectively avoiding numerical blow-up. In particular, the sG projection coefficients associated with $\tilde{D}(z)$ remain finite even when the corresponding integrals involving the original singular coefficient $D(z)$ diverge, making the sG approximation well-defined and its numerical simulation more stable. The function $\Psi(r)$ is bounded and smooth, as will be proved in Section 4.1. This strategy is the core of the str-RB algorithm and plays an essential role in guaranteeing the numerical stability of the algorithm. After regularizing the diffusion term as $\tilde{D}$, the particle system Eq. \eqref{2} becomes
\begin{equation}
\left\{
\begin{aligned}
dX_i(z, t) &= V_i(z, t)\, dt,\\
dV_i(z, t) &= \frac{1}{N} \sum_{j=1}^N P(z, X_i, X_j)\,(V_j - V_i)\, dt 
              + \sqrt{2\tilde{D}(z)}\, dW_i,
\end{aligned}
\right.
\qquad i=1,\dots, N,
\label{3}
\end{equation}
which allows the sG projection to be performed stably.
\subsection{sG approximation}
\par With the regularized diffusion term $\tilde{D}$, the sG approximation of Eq. \eqref{3} is given by the gPC expansion \cite{3}
\begin{align}
X_i^M(z,t) = \sum_{m=0}^M \hat{X}_{i,m}(t)\Phi_m(z), \quad
V_i^M(z,t) = \sum_{m=0}^M \hat{V}_{i,m}(t)\Phi_m(z), \quad i = 1,\dots,N,\label{6}
\end{align}
where $\{\Phi_m(z)\}_{m=0}^M$ is a set of orthogonal polynomials with respect to the probability density function $p(z)$. Here, $\hat{X}_{i,m}$ and $\hat{V}_{i,m}$ are the projections of the solutions onto $\Phi_m$, namely,
\begin{align*}
\hat{X}_{i,m}(t) = \mathbb{E}\bigl[X_i(z,t)\Phi_m(z)\bigr], \quad
\hat{V}_{i,m}(t) = \mathbb{E}\bigl[V_i(z,t)\Phi_m(z)\bigr].
\end{align*}
Inserting \eqref{6} into \eqref{3} and applying a Galerkin projection yields the particle sG approximation of Eq. \eqref{3}
\begin{align}
\begin{cases}
d\hat{X}_{i,h}(t) = \hat{V}_{i,h}(t) dt, \\
d\hat{V}_{i,h}(t) = \displaystyle\frac{1}{N}\sum_{j=1}^N \sum_{k=0}^M P_{hk}^{ij}(\hat{V}_{j,k} - \hat{V}_{i,k}) dt + \tilde{D}_h dW_i,\label{4}
\end{cases}
\end{align}
with $P_{hk}^{ij} = \mathbb{E}\bigl[P(z,X_i^M,X_j^M)\Phi_h(z)\Phi_k(z)\bigr]$ and $\tilde{D}_h = \mathbb{E}\bigl[\sqrt{2\tilde{D}(z)}\Phi_h(z)\bigr]$. This is a $2N(M+1)$-dimensional SDE with deterministic coefficients.
\subsection{Random batch approach}
\par Solving Eq. \eqref{4} directly involves $\mathcal{O}(N^2)$ particle interactions and is therefore computationally expensive. To reduce the cost, a variant of the RBM \cite{22} is adopted to solve the deterministic-coefficient SDE Eq. \eqref{4}. At each time step, an $S$-element subset $S_{i}$ is randomly selected from the index set $\{1,\dots,N\}$ for each particle $i$, and particle $i$ interacts only with $S$ particles in $S_{i}$ during that step. After applying the random batch approach, the SDE \eqref{4} reads
\begin{align}
\begin{cases}
d\hat{X}_{i,h}(t) = \hat{V}_{i,h}(t) dt, \\
d\hat{V}_{i,h}(t) = \displaystyle\frac{1}{S}\sum_{j\in S_{i}} \sum_{k=0}^M P_{hk}^{ij}(\hat{V}_{j,k} - \hat{V}_{i,k}) dt + \tilde{D}_h dW_i,\label{14}
\end{cases}    
\end{align}
For convenience, we denote 
\begin{align}
F_{i,h}:=\displaystyle\frac{1}{S}\sum_{j\in S_{i}} \sum_{k=0}^M P_{hk}^{ij}(\hat{V}_{j,k} - \hat{V}_{i,k}) .
\end{align} 
The random batch approach reduces the per-step computational complexity significantly.
\subsection{Modified tamed scheme}
\par After the computational cost has been reduced by the random batch approach, the modified tamed scheme \cite{2} is employed here to solve Eq. \eqref{14}, since $P(z, \cdot)$ in Eq. \eqref{2} may be non-smooth or singular. Specifically, for better stability, we introduce the tamed drift term
\begin{align*}
    \tilde{F}_{i,h}=\mathcal{T}_{\varphi}(F_{i,h})
\end{align*}
where the taming function $\mathcal{T}_{\varphi}(r)$ is defined as
\begin{align}
\mathcal{T}_{\varphi}(r)=\frac{r}{1 + \varphi\bigl(\gamma_{\text{drift}} \tau^{\alpha_{\text{drift}}} |r|\bigr)},\label{11}
\end{align}
with $\tau$ denoting the time step, $\alpha_{\text{drift}} \in (0, \frac{1}{2}]$, $\gamma_{\text{drift}} \geq 0$, and $\varphi: \mathbb{R}_+ \to \mathbb{R}_+$ a cut-off function satisfying:
\begin{align}    
\varphi(\ell) = 
\begin{cases}
0, & \ell \leq 1, \\
\ell, & \ell \geq 2, \\
\text{smooth increasing}, & \ell \in (1, 2). 
\end{cases}\label{12*}
\end{align}
The SDE with deterministic coefficients then becomes
\begin{align}
\begin{cases}
d\hat{X}_{i,h}(t) = \hat{V}_{i,h}(t) dt, \\
d\hat{V}_{i,h}(t) = \tilde{F}_{i,h} dt + \tilde{D}_h dW_i,\label{15}
\end{cases}    
\end{align}
which is finally discretized by the Euler-Maruyama scheme. The taming strategy mitigates the numerical instability caused by the singular drift term. 
\subsection{Algorithm}
\par Now the complete algorithm is presented. First and most critically, as described in Section 3.1, the exponential regularization strategy is adopted to the singular diffusion term before the sG projection , ensuring that the projection can be performed stably. The sG projection then converts the random-coefficient system into a SDE with deterministic coefficients, after which a random batch approach is employed to reduce computational complexity, and the modified tamed scheme is used to solve the resulting SDE, mitigating the instability induced by the drift term. This sequential treatment is essential because the singular diffusion term must be regularized first to guarantee a stable sG projection, the RBM then reduces the computational complexity, and the tamed scheme finally addresses the singularities in the drift term. 
\par The complete str-RB algorithm is given in Algorithm \ref{Algorithm2}.
\begin{FloatAlgo}{stochastic tamed-regularized Random Batch (str-RB) algorithm}
\label{Algorithm2}
\begin{algorithmic}[1]
\STATE Consider $N$ samples $(X_i, V_i)$ from $f_0(x, v)$ and fix $S \leq N$.
\STATE Perform gPC representation on the particles: $(\hat{X}_{i,h}, \hat{V}_{i,h})$, for $h = 0, \dots, M$.
\FOR{$n = 0, \dots, \lfloor T/\tau \rfloor-1$}
    \STATE Generate $N$ random variables $\{\eta_i\}_{i=1}^N \sim \mathcal{N}(0, 1)$.
    \FOR{$i = 1, \dots, N$}
        \STATE Sample $S$ particles $\{j_1, \dots, j_S\} := S_i$ uniformly without repetition.
        \STATE Update the position and velocity:
        \begin{align*}
        \hat{X}_{i,h}^{n+1} &= \hat{X}_{i,h}^n + \hat{V}_{i,h}^n \tau, \\
        \hat{V}_{i,h}^{n+1} &= \hat{V}_{i,h}^n + \tau \tilde{F}_{i,h} + \sqrt{\tau} \tilde{D}_h \eta_i.
        \end{align*}
    \ENDFOR
\ENDFOR
\end{algorithmic}
\end{FloatAlgo}
Introducing the exponential regularization strategy is the core step of the str-RB algorithm, which is designed to mitigate numerical blow-up induced by singularities in the diffusion term and outperforms taming-type treatments for  highly singular diffusion terms. To be specific, recall from Eqs. \eqref{3-4} and \eqref{11} that, for sufficiently large $r$, the attenuation factors of the exponential regularization and of the taming scheme are respectively
\begin{align*}
\psi(r)=\exp(-\tau^{\alpha_{\text{diff}}}\gamma_{\text{diff}} r^{\beta_{\text{diff}}}), \qquad h(r)=\frac{\mathcal{T}_{\varphi}(r)}{r}=\frac{1}{1+\tau^{\alpha_{\text{drift}}}\gamma_{\text{drift}} r}.
\end{align*}
We have two observations when the powers of $\tau$ coincide, i.e., $\alpha_{\text{drift}}=\alpha_{\text{diff}}=\alpha$, with $\gamma_{\text{drift}}, \gamma_{\text{diff}}>0$. First, to address the difficulties caused by Brownian motion, the diffusion term requires a stronger treatment than the taming strategy. Since $\psi(r)$ decays exponentially, one has $\psi(r)<h(r)$ for sufficiently large $r$ , and thus $\psi(r)$ provides stronger suppression and ensures better stability. Second, as $\tau\to 0$, one has 
\begin{align*}
\psi(r)\sim 1-\tau^\alpha\gamma_{\text{diff}} r^{\beta_{\text{diff}}},\quad  h(r)\sim 1-\tau^\alpha\gamma_{\text{drift}}r , 
\end{align*}
and $r^{\beta_{\text{diff}}}\ll r$ for large $r$, so the exponential regularization restores the unperturbed behavior more accurately as $\tau \to 0$,  yielding a closer approximation to the solution of the original system. In view of these two observations, exponential regularization is a natural choice for handling numerical instability induced by strong singularities in the diffusion term in the presence of Brownian motion.
\par Moreover, for smooth cases, setting $ \gamma_{\text{drift}} = \gamma_{\text{diff}}= 0 $ reduces the str-RB algorithm to the standard sG algorithm, i.e., Algorithm 3 in \cite{3}. The str-RB algorithm is therefore consistent with the method of \cite{3} and retains its ability to solve smooth problems. The effectiveness and accuracy of the proposed algorithm will be demonstrated through numerical tests in Section 5.
\section{Theoretical properties}
\par This section focuses on the theoretical properties of the str-RB algorithm. Specifically, the numerical stability, computational complexity, and regularization-induced model perturbation are analyzed in detail, followed by the discussion of the numerical error in the regularized system.
\subsection{Numerical stability}
\par We begin with the numerical stability analysis. The exponential regularization strategy plays a critical role in reducing errors and maintaining numerical stability. First we establish the smoothness of $\Psi(r)$ in Eq. \eqref{3-3} rigorously, which serves as an important theoretical guarantee of the regularization strategy.
\begin{proposition}
    Both $\Psi(r)$ and its derivative are bounded.
\end{proposition}
\begin{proof}
\par We split the proof into two parts. We first establish the boundedness of $\Psi(r)$ and then prove the boundedness of $\Psi'(r)$.
\par \textbf{Boundedness of $\Psi(r)$.} Recall from Eq. \eqref{3-3} that
\begin{align}
\Psi(r)=r\exp\left( -A  r^{\beta_{\text{diff}}} \right),\label{4-1}    
\end{align}
where we denote $A:=\tau^{\alpha_{\text{diff}}} \gamma_{\text{diff}}>0$. 
\par A direct computation yields
\begin{align}
    \Psi'(r) = \exp\left(-A r^{\beta_{\text{diff}}}\right) \cdot \left(1 - A\beta_{\text{diff}} r^{\beta_{\text{diff}}}\right). \label{10}
\end{align}
It follows that 
\begin{align*}
    \Psi'(r)=0 \Rightarrow r=\left( \frac{1}{A \beta_{\text{diff}}} \right)^{1/\beta_{\text{diff}}}:=r_0.
\end{align*}
\par Since $\Psi'(r)>0$ for $0\leq r<r_0$ and $\Psi'(r)<0$ for $r>r_0$, $\Psi(r)$ attains its maximum at $r_0$. Thus, $0 \le \Psi(r) \le \Psi(r_0)$ for all $r \ge 0$.
\par \textbf{Boundedness of $\Psi'(r)$.} Let $r_1=A r^{\beta_{\text{diff}}}$. Substituting into Eq. \eqref{10}, one gains $\Psi'(r)=h(r_1)$, where
\begin{align*}
   h(r_1)=e^{-r_1} (1-\beta_{\text{diff}} r_1).
\end{align*}
\par For all $x>0$, $0<e^{-x}<1$, and $0<xe^{-x}\le e^{-1}$, which gives 
\begin{align*}
    |h(r_1)|\le|e^{-r_1}|+\beta_{\text{diff}}|r_1 e^{-r_1}|\le 1+\beta_{\text{diff}} e^{-1},
\end{align*}
indicating the boundedness of $\Psi'(r)$. 
\par Combining the two parts, one can conclude that $\Psi(r)$ and its derivative are bounded.
\end{proof}
\par The smoothness of $\Psi(r)$ is crucial for the stability of the algorithm. The boundedness of $\Psi(r)$ implies the boundedness of the regularized diffusion term, which prevents large fluctuations driven by large values of $D(z)$ and Brownian motion and thus avoids numerical blow-up. The boundedness of  $\Psi'(r)$ ensures that $\Psi$ does not exhibit significant fluctuations, which effectively mitigates the impact of sampling errors in the approximation of the expectation $\tilde{D}_h = \mathbb{E}\biggl[\sqrt{2\tilde{D}(z)}\Phi_h(z)\biggr]$.  Moreover, $\Psi(r) \to r~\text{as the time step}~\tau \to 0$, so the exponential regularization does not introduce excessive perturbation to the original system for small time steps. Consequently, the exponential regularization strategy serves as a key stabilization safeguard of the str-RB algorithm, enabling it to handle highly singular diffusion terms that are difficult for existing methods.
\par Next, the contribution of the modified tamed scheme to the numerical stability of the str-RB algorithm is analyzed. For the SDE obtained after sG projection and random batch approximation, a direct numerical solution may still encounter numerical instability due to potential singularities of $P(z,\cdot)$ in Eq. \eqref{2}. As shown in Section 3.4, to overcome numerical instability, the modified tamed scheme is employed to solve the SDE Eq. \eqref{14}, complementing the exponential regularization strategy applied to the diffusion term. We now establish the smoothness of the taming function $\mathcal{T}_{\varphi}(r)$.
\begin{proposition}
    If the cut-off function $\varphi(\ell)$ satisfies $|\varphi'(\ell)|\le C_0$ for some constant $C_0$, then both the taming function $\mathcal{T}_{\varphi}(r)$ and its derivative are bounded.
\end{proposition}
\begin{proof}
\par Since $\mathcal{T}_{\varphi}(-r) = -\mathcal{T}_{\varphi}(r)$, it suffices to consider $r\ge 0$. We write
\begin{align*}
    \mathcal{T}_{\varphi}(r) = \dfrac{r}{1+\varphi(cr)},\quad c=\gamma_{\text{drift}}\tau^{\alpha_{\text{drift}}}.
\end{align*}
The proof proceeds in two steps. We first show that the taming function $\mathcal{T}_{\varphi}(r)$ is bounded and subsequently establish the boundedness of its derivative.
\par \textbf{Boundedness of $\mathcal{T}_{\varphi}(r)$.} 
We consider three cases according to the range of $r$.
\par (1) If $0\le r\le 1/c$, then $cr\le 1$, so the cut-off function 
 $\varphi(cr)=0$ and thus $\mathcal{T}_{\varphi}(r)=r\in[0,1/c]$.
\par (2) If $r\ge 2/c$, then $cr\ge 2$, hence $\varphi(cr)=cr$. Then one yields that $\mathcal{T}_{\varphi}(r)=\dfrac{r}{1+cr}\in(0,1/c)$.
\par (3) If $1/c<r<2/c$, $\varphi(cr)>0$ since $\varphi$ is increasing on $(1,2)$. Thus one has 
\begin{align*}
0< \mathcal{T}_{\varphi}(r)<r<2/c.  
\end{align*}
\par Combining the three cases above, $\mathcal{T}_{\varphi}(r)$ is bounded.
\par \textbf{Boundedness of the derivative.}
Direct computation gives
\begin{align*}    
\bigl(\mathcal{T}_{\varphi}\bigr)'(r) = \dfrac{1+\varphi(\ell)-\ell\,\varphi'(\ell)}{(1+\varphi(\ell))^2}:=B(\ell),
\end{align*}
where $\ell=cr$. Similar to the previous step, three cases are considered as follows.
\par (1) If $\ell\le 1$, then $\varphi(\ell)=0$, $\varphi'(\ell)=0$, so $B(\ell)=1$.
\par (2) If $\ell\ge 2$, one has $\varphi(\ell)=\ell$, $\varphi'(\ell)=1$, which gives $B(\ell)=\dfrac{1}{(1+\ell)^2}\in(0,1]$.
\par (3) If $1<\ell<2$, as $|\varphi'(\ell)|\le C_0$, one has $B_{-}(\ell)\le B(\ell)\le B_{+}(\ell)$,
where
\begin{align}
    B_{\pm}(\ell)=\dfrac{1+\varphi(\ell)\pm 2C_0}{(1+\varphi(\ell))^2}.\label{11*}
\end{align}
As can be seen from Eq. \eqref{11*}, both $B_{-}(\ell)$ and $B_{+}(\ell)$ are continuous on $[1,2]$, which implies that they are bounded on $[1,2]$. Therefore, $B(\ell)$ is bounded on $(1,2)$.
\par In summary,  $\bigl(\mathcal{T}_{\varphi}\bigr)'(r)=B(\ell)$ is bounded.
\end{proof}
\par The smoothness of the taming function ensures that the tamed SDE with deterministic coefficients, i.e., Eq. \eqref{15}, can be solved stably using the Euler-Maruyama scheme. Moreover,  $\mathcal{T}_{\varphi}(r)\to r ~\text{as}~\tau \to 0$, so the modified tamed scheme does not introduce excessive perturbations to the original system for small time steps.
\par It is noted in \cite{2} that the modified tamed scheme is designed for SDEs with non-globally Lipschitz drift terms. In this work, the drift term of the SDE obtained after sG projection may in fact be globally Lipschitz but with a very large Lipschitz norm. The modified tamed scheme is still employed rather than solving the SDE directly to further mitigate the effect of the non-smooth drift term in Eq. \eqref{2} and to prevent the amplification of errors, such as those arising from the computation of the expectation $P_{hk}^{ij}$ via numerical quadrature or Monte Carlo methods. 
\subsection{Computational complexity}
\par Having established the stability of the str-RB algorithm, we now analyze its computational complexity. In the str-RB algorithm, a variant of the RBM is adopted, where each particle interacts only with a small random subset of the other particles. The reduction in computational complexity is quantified in the following proposition.
\begin{proposition}
The random batch approach in Eq. \eqref{14} reduces the per-step computational complexity of solving the SDE obtained after the sG projection from $\mathcal{O}(M^2N^2)$ to $\mathcal{O}(M^2SN)$.
\end{proposition}
\begin{proof}
Following the complexity analysis in \cite{3}, a direct solve of the SDE after sG projection requires each particle to interact with all other $(N-1)$ particles, giving  $\mathcal{O}(N^{2})$ interactions per step.  Since each particle has $(M+1)$ gPC coefficients, each interaction costs $\mathcal{O}(M^2)$. The per-step complexity of a direct solve of the SDE is therefore $\mathcal{O}(M^2N^2)$.
\par In the str-RB algorithm, the random batch approach is adopted. Particle $i$ interacts only with $S$ particles in $S_{i}$, which reduces the number of particle interactions per step to $\mathcal{O}(SN)$.
As the cost per interaction remains  $\mathcal{O}(M^2)$, the per-step complexity of the str-RB algorithm is 
\begin{align*}
    \mathcal{O}(M^2SN).
\end{align*}
\end{proof}
\par The per-step complexity $\mathcal{O}(M^2SN)$ of the str-RB algorithm corresponds to that in \cite{3}. For example, taking $N = 1000$, $M = 2$, and $S = 5$, consistent with the numerical tests in Section 5, the per-step complexity is reduced from $\mathcal{O}\left(M^2 N^2\right) = \mathcal{O}(4\times 10^{6})$  to $\mathcal{O}\left(M^2 S N\right) = \mathcal{O}(2\times 10^{4})$, i.e., by a factor of $200$. 
\par As noted in \cite{3}, applying the random batch approach to the deterministic-coefficient SDE obtained after sG projection introduces an additional error of $O\left(\sqrt{\frac{1}{S}-\frac{1}{N}}\right)$ compared to a direct solve of the SDE. However, when the coefficients in the drift and diffusion terms are singular with respect to the random variable, numerical instability induced by singular coefficients becomes the dominant source of error, exceeding the random batch error by several orders of magnitude. Therefore, the additional error introduced by the random batch strategy does not significantly affect the performance of the str-RB algorithm, while the reduction in computational complexity is substantial.
\subsection{Analysis of the regularization-induced model perturbation}
\par In this section, we analyze the asymptotic behavior, as the time step tends to zero,  of  the expectation and variance of the regularized diffusion coefficient, which characterizes the perturbation introduced by the regularization. We first state the assumptions on the random variable and the diffusion term.
\begin{assumption}
\label{assumption1}

\par (1) The random variable $z$ follows the uniform distribution on $(-1,1)$. 
\par (2) The singular diffusion coefficient has the form
\begin{align}\label{eq:Dz}
D(z)=a(z)(1-z)^{-K},
\end{align}
where $a(z)>0,~a(z)\in C^{\infty}$ and $K$ is a positive integer. 
\end{assumption}
\par Recall that the regularized diffusion coefficient is defined as
\begin{align*}
\tilde{D}(z)=D(z)\psi(D(z)),
\end{align*}
where
\begin{align*}
\psi(r)=\exp(-Ar^{\beta_{\text{diff}}}),
\end{align*}
with $A:=\tau^{\alpha_{\text{diff}}}\gamma_{\text{diff}}>0$. The following proposition describes the asymptotic behavior of the expectation and variance of $\tilde{D}$ as $\tau\rightarrow0$.
\begin{proposition}
\label{proposition:model_perturbation}
Under Assumption \ref{assumption1}, with $\alpha_{\text{diff}},\gamma_{\text{diff}},\beta_{\text{diff}}>0$, as $\tau\rightarrow0$, the expectation and variance of the regularized diffusion term $\tilde{D}$ satisfy
\begin{align}
\mathbb{E}[\tilde{D}]
=
\begin{cases}
O\left(
\tau^{-\frac{\alpha_{\text{diff}}(K-1)}{K\beta_{\text{diff}}}}
\right),
&K>1,\\
O\left(\log\frac{1}{\tau}\right),
&K=1,
\end{cases}\label{17}
\end{align}
and
\begin{align}
\mathrm{Var}(\tilde{D})
=
O\left(
\tau^{-\frac{\alpha_{\text{diff}}(2K-1)}{K\beta_{\text{diff}}}}
\right).\label{18}
\end{align}
\end{proposition}
\begin{proof}
\par Since $a(z) \in C^{\infty}$, it suffices to prove it for the case $a(z)=1$, i.e., $D(z)=(1-z)^{-K}$. The proof consists of two parts. We first analyze the asymptotic behavior of the expectation Eq. \eqref{17}, and then study the second moment, from which the variance estimate Eq. \eqref{18} follows.

\par Since $z\sim U(-1,1)$, the expectation of $\tilde{D}$ is given by
\begin{align}
\mathbb{E}[\tilde{D}]
=
\frac{1}{2}
\int_{-1}^{1}
(1-z)^{-K}
\exp\left(-A(1-z)^{-K\beta_{\text{diff}}}\right)dz ,\label{19}
\end{align}
where $A:=\tau^{\alpha_{\text{diff}}}\gamma_{\text{diff}}$, hence $\tau \to 0$ is equivalent to $A \to 0$. Introducing the change of variables $y=A(1-z)^{-K\beta_{\text{diff}}}$, one obtains
\begin{align}
1-z=
\left(\frac{A}{y}\right)^{\frac{1}{K\beta_{\text{diff}}}},\label{20}
\end{align}
and
\begin{align}
dz
=
\frac{1}{K\beta_{\text{diff}}}
A^{\frac{1}{K\beta_{\text{diff}}}}
y^{-\frac{1}{K\beta_{\text{diff}}}-1}dy .\label{21}
\end{align}
Substituting Eqs. \eqref{20} and \eqref{21} into Eq. \eqref{19} gives
\begin{align}
\mathbb{E}[\tilde{D}]
&=
\frac{1}{2K\beta_{\text{diff}}}
A^{-\frac{K-1}{K\beta_{\text{diff}}}}
\int_{2^{-K\beta_{\text{diff}}}A}^{+\infty}
y^{\frac{K-1}{K\beta_{\text{diff}}}-1}e^{-y}dy. \label{22*}
\end{align}
\textbf{Case 1:} For $K>1$, the integral term in Eq. \eqref{22*} can be bounded by the Gamma function
\begin{align*}
\int_{2^{-K\beta_{\text{diff}}}A}^{+\infty}
y^{\frac{K-1}{K\beta_{\text{diff}}}-1}e^{-y}dy
\leq
\Gamma\left(\frac{K-1}{K\beta_{\text{diff}}}\right),
\end{align*}
where $\Gamma(\cdot)$ denotes the Gamma function. Hence,
\begin{align}
\mathbb{E}[\tilde{D}]
=
O\left(
A^{-\frac{K-1}{K\beta_{\text{diff}}}}
\right).\label{22}
\end{align}
\textbf{Case 2:} For $K=1$, Eq. \eqref{22*} becomes 
\begin{align*}
\mathbb{E}[\tilde{D}]=
\frac{1}{2\beta_{\text{diff}}} \int_{2^{-\beta_{\text{diff}}}A}^{+\infty}y^{-1}e^{-y}dy.   
\end{align*}
Using the asymptotic expansion of the exponential integral function
\begin{align*}
\int_{x}^{+\infty}y^{-1}e^{-y}dy
=
-\eta -\log x+O(x),
\qquad x\rightarrow0^{+},
\end{align*}
where $\eta$ is the Euler constant, one derives that as $A\to 0$,
\begin{align}
\mathbb{E}[\tilde{D}]
=
O\left(
\log\frac{1}{A}
\right).\label{22n}
\end{align} 
Thus, we obtain the asymptotic behavior of the expectation of the
regularized diffusion coefficient.
\par Next, we consider the variance. Similarly, 
\begin{align*}
\mathbb{E}[\tilde{D}^{2}]
=
\frac{1}{2}
\int_{-1}^{1}
(1-z)^{-2K}
\exp\left(-2A(1-z)^{-K\beta_{\text{diff}}}\right)dz .
\end{align*}
A transformation $y=2A(1-z)^{-K\beta_{\text{diff}}}$ gives
\begin{align*}
\mathbb{E}[\tilde{D}^{2}]
&=
\frac{1}{2K\beta_{\text{diff}}}
(2A)^{-\frac{2K-1}{K\beta_{\text{diff}}}}
\int_{2^{1-K\beta_{\text{diff}}}A}^{+\infty}
y^{\frac{2K-1}{K\beta_{\text{diff}}}-1}e^{-y}dy .
\end{align*}

Since the integral term is bounded by $\Gamma\left(\frac{2K-1}{K\beta_{\text{diff}}}\right)$, it follows that
\begin{align}
\mathbb{E}[\tilde{D}^{2}]
=
O\left(
A^{-\frac{2K-1}{K\beta_{\text{diff}}}}
\right).\label{23}
\end{align}
Combining Eqs. \eqref{22}, \eqref{22n} and \eqref{23}, we have
\begin{align}
\mathrm{Var}(\tilde{D})
&=
\mathbb{E}[\tilde{D}^{2}]
-
(\mathbb{E}[\tilde{D}])^{2}
\nonumber\\
&=
O\left(
A^{-\frac{2K-1}{K\beta_{\text{diff}}}}
\right).\label{25}
\end{align}
Finally, substituting $A=\tau^{\alpha_{\text{diff}}}\gamma_{\text{diff}}$ into Eqs. \eqref{22}, \eqref{22n}, and \eqref{25},  and noting that negative powers of $\gamma_{\mathrm{diff}}$ and the term $-\log\gamma_{\mathrm{diff}}$ are constant with respect to $\tau$, we conclude that, as $\tau \to 0$, for $K>1$,
\begin{align}
\mathbb{E}[\tilde{D}]
=  O\left(\gamma_{\text{diff}}^{-\frac{K-1}{K\beta_{\text{diff}}}}
\tau^{-\frac{\alpha_{\text{diff}}(K-1)}{K\beta_{\text{diff}}}}
\right)=O\left(
\tau^{-\frac{\alpha_{\text{diff}}(K-1)}{K\beta_{\text{diff}}}}
\right);  \label{26}
\end{align}
 for $K=1$,
\begin{align}
\mathbb{E}[\tilde{D}]
=  O\left(\log\frac{1}{\tau^{\alpha_{\text{diff}}}\gamma_{\text{diff}}}\right)=O\left(\alpha_{\text{diff}}\log\frac{1}{\tau}-\log \gamma_{\text{diff}}\right)=O\left(\log\frac{1}{\tau}\right).\label{27}
\end{align}
And
\begin{align}
\mathrm{Var}(\tilde{D})
=
O\left(\gamma_{\text{diff}}^{-\frac{2K-1}{K\beta_{\text{diff}}}}
\tau^{-\frac{\alpha_{\text{diff}}(2K-1)}{K\beta_{\text{diff}}}}
\right)=O\left(
\tau^{-\frac{\alpha_{\text{diff}}(2K-1)}{K\beta_{\text{diff}}}}
\right). \label{28}
\end{align}
Combining Eqs. \eqref{26}, \eqref{27}, and \eqref{28} completes the proof.
\end{proof}
\par Proposition \ref{proposition:model_perturbation} characterizes the perturbation introduced by the exponential regularization through the asymptotic behavior of the moments. It shows that the regularized coefficient progressively recovers the statistical behavior of the original diffusion coefficient as $\tau \to 0$. The choice of the regularization parameter therefore involves a trade-off: a weaker regularization reduces the modeling error but may weaken the stabilization effect, whereas a stronger regularization enhances numerical stability at the expense of a larger model perturbation.
\subsection{Discussion of numerical errors for the regularized system}
\par For the original system Eq. \eqref{2}, our algorithm consists of two stages: regularization and numerical solution of the regularized system. The regularization-induced modeling error arising from the first stage is analyzed in Section 4.3, while this section discusses the numerical errors in the second stage. For the convenience of the discussion, we assume in this section that the drift term $P$ is smooth and only the diffusion term is singular. Specifically, assume that $P(z,\cdot)$ satisfies
\begin{align}
P(z, \cdot) \in C^{\infty},    \label{4-18}
\end{align}
and for any integer $m\geq 0$,
\begin{align}
\int_{-1}^{1}
\left|
\frac{\partial^{m}P(z, \cdot)}{\partial z^{m}}
\right|^{2}
dz
\leq M_m ,\label{4-19}
\end{align}
where $M_m$ is a constant depending only on $m$.
\par We now discuss the main component‌s of numerical errors in solving the regularized system Eq. \eqref{3}, including sG projection error, error of the random batch approach, and error of discretization. 
\par (1) sG projection error.
\par Following the notation in Eq. \eqref{4-1}, and noting that every derivative of the exponential function retains the exponential factor, for any integer $j\geq 0$, one has
\begin{align}
\Psi^{(j)}(r)
=
Q_j(r)\exp(-Ar^{\beta_{\text{diff}}}),\label{12}
\end{align}
where $Q_j$ is a generalized polynomial of the form
\begin{align*}
Q_j(r) = \sum_{k=1}^{m_j} c_{jk} r^{\eta_{jk}},\quad m_j\in\mathbb{N}^*,\ c_{jk}\in\mathbb{R},\ \eta_{jk}\in\mathbb{R}.
\end{align*}
On the other hand, it follows from \eqref{eq:Dz} that
\begin{align}
D^{(j)}(z)
&=
\sum_{\ell=0}^{j}
\binom{j}{\ell}
a^{(\ell)}(z)(K)_{j-\ell}
(1-z)^{-K-j+\ell}
\nonumber\\
&=
B_j(z)(1-z)^{-K-j},\label{13}
\end{align}
where
\begin{align*}
B_j(z)=
\sum_{\ell=0}^{j}
\binom{j}{\ell}
a^{(\ell)}(z)(K)_{j-\ell}(1-z)^\ell,
\end{align*}
and $(K)_n=K(K+1)\cdots(K+n-1)$ for $n\geq1$, with $(K)_0=1$. Since $a(z)\in C^\infty$, $B_j(z)$ is smooth and bounded on $[-1,1]$.

Combining Eqs. \eqref{12} and \eqref{13} and applying the Fa\`a di Bruno formula to the composite function $\tilde{D}(z)=\Psi(D(z))$, one obtains
\begin{align}
\frac{d^s\tilde{D}(z)}{dz^s}
=
\Pi_s(z)
\exp\left(-Aa(z)^{\beta_{\text{diff}}}(1-z)^{-K\beta_{\text{diff}}}\right),\label{4.20}
\end{align}
where 
\begin{align*}
\Pi_s(z)= \sum_{k=1}^{n_s} b_{sk}(z) (1-z)^{-\zeta_{sk}},\quad n_s\in\mathbb{N}^*,\ b_{sk}(z)\in C^{\infty},\ \zeta_{sk}\geq0.   
\end{align*}
Since $a(z)>0$ and $a(z)\in C^\infty$, the exponential decay in Eq. \eqref{4.20} dominates any algebraic growth, and it follows that
\begin{align*}
\frac{d^s\tilde{D}(z)}{dz^s}\in L^2(-1,1)
\end{align*}
for an arbitrary integer $s\geq0$. Therefore,
\begin{align*}
\tilde{D}(z)\in C^\infty(-1,1)\cap H^{s}(-1,1),
\quad \forall~s\in\mathbb{N}_{+}.
\end{align*}
Together with the regularity assumptions Eqs. \eqref{4-18} and \eqref{4-19} on $P(z,\cdot)$, both the drift and diffusion coefficients of the regularized system possess square-integrable derivatives of arbitrary order with respect to the random variable $z$. According to the standard stochastic Galerkin approximation theory \cite{31}, the sG projection error of Eq. \eqref{3} exhibits spectral accuracy of order $O(M^{-J})$, where $M$ is the gPC truncation order in Eq. \eqref{6} and $J$ is an arbitrary positive integer.

(2) Error of the random batch approach. According to the result in \cite{3}, the error introduced by the random batch approach when solving Eq. \eqref{4} has the order $O\left(\sqrt{\frac{1}{S}-\frac{1}{N}}\right)$, which vanishes when $S=N$.

(3) Error of discretization. Since both the drift and diffusion terms of the regularized system are $C^{\infty}$, the discretization error of solving Eq. \eqref{15} has the same order as that of the Euler--Maruyama scheme.
\par The analysis in Section 4.3 and the above discussion together imply that the convergence behavior observed in the numerical experiments of Section 5 reflects the combined influence of the numerical error of the regularized system and the regularization-induced model perturbation. Due to the trade-off between these two error components, the convergence behavior of the original singular system in the experiments differs from that of the regularized system alone.
\section{Numerical experiments and results analysis}
\par In this section, a series of numerical tests is carried out to validate the convergence, accuracy, and numerical stability of the str-RB algorithm, i.e.,  Algorithm \ref{Algorithm2}. We primarily consider cases where the drift term $P(z)$ and the diffusion term $D(z)$ have singularities, particularly those where the diffusion term exhibits challenging singularities that existing methods struggle to handle.
\subsection{Experimental setup}
Before presenting detailed numerical results, we describe the benchmark algorithms, basic parameters and error metrics adopted in the tests. 
\par The following four algorithms are compared to validate the effectiveness of the str-RB algorithm.
\begin{itemize}
    \item \textbf{s-RB}: the original baseline, i.e., the particle sG algorithm proposed in \cite{3} (Algorithm 3 in Section 7.1 therein), without any additional treatment of $P(z,\cdot)$ and $D(z)$;
    \item \textbf{st-RB}: an extension of s-RB by applying the modified tamed scheme in \cite{2} to solve the SDE with deterministic coefficients obtained after the sG projection;
    \item \textbf{stt-RB}: an enhanced baseline constructed in this work, which extends st-RB by additionally taming the diffusion term $D(z)$,  that is, replacing $D$ with the tamed version $\mathcal{T}_{\varphi_\text{diff}}(D)$ before the sG projection. See Table \ref{para} for the definition of $\varphi_\text{diff}$;
    \item \textbf{str-RB}: the stochastic tamed-regularized Random Batch algorithm proposed in this work, i.e., Algorithm \ref{Algorithm2}.
\end{itemize}
\par In all tests, we set $z\sim\mathcal{U}(-1,1)$ and adopt Legendre polynomials as the basis functions for the gPC expansion. Unless otherwise noted, the following parameters are used: the number of particles $N=1000$, gPC order $M=2$ in Eq. \eqref{6}, batch size $S=5$, final time $T=1$, and time step $\tau=0.004$. 
\par The algorithm-specific parameters and cut-off functions  are summarized in Table \ref{para}. Unless explicitly stated otherwise, all subsequent numerical tests are conducted with the parameter values given in this table. It is noted that the specific choices of the parameter values have a relatively significant influence on the performance of the algorithm, whereas the precise form of the cut-off function on $(1,2)$ has a comparatively minor influence.
\begin{table}[htbp]
  \centering
  \caption{Parameters for different algorithms}
  \begin{tabular}{@{}ll@{}}
    \toprule
    Algorithm & Parameters \\
    \midrule
    st-RB & 
    \begin{tabular}[c]{@{}l@{}}
     Modified tamed scheme Eq. \eqref{11} for SDE obtained after sG projection: \\
    $\alpha_{\text{drift}} = 0.25$, $\gamma_{\text{drift}} = 2$, \\
    cut-off function 
    $
    \varphi_{\text{drift}}(\ell)=
    \begin{cases}
     0, & \ell \le 1, \\
     \ell, & \ell \ge 2, \\
     \ell^3 - 2\ell^2 + \ell, & \ell \in (1,2).
     \end{cases}$
     \end{tabular} \\
    \midrule
    stt-RB & 
    \begin{tabular}[c]{@{}l@{}}
      Modified tamed scheme Eq. \eqref{11} for SDE obtained after sG projection:\\
      $\alpha_{\text{drift}} = 0.25$, $\gamma_{\text{drift}} = 2$, \\
      cut-off function
      $
      \varphi_{\text{drift}}(\ell)=
    \begin{cases}
     0, & \ell \le 1, \\
     \ell, & \ell \ge 2, \\
     \ell^3 - 2\ell^2 + \ell, & \ell \in (1,2).
     \end{cases}
      $\\
      Modified tamed scheme Eq. \eqref{11} for the diffusion term: \\
      $\alpha_{\text{diff}} = 0.25$, $\gamma_{\text{diff}} = 1$, \\
      cut-off function 
      $
      \varphi_{\text{diff}}(\ell)=
      \begin{cases}
        0, & \ell \le 1, \\
        \ell, & \ell \ge 2, \\
        2\ell-2, & \ell \in (1,2).
      \end{cases}
      $
    \end{tabular} \\
    \midrule
    str-RB & 
    \begin{tabular}[c]{@{}l@{}}
    Modified tamed scheme Eq. \eqref{11} for SDE obtained after sG projection:\\
      $\alpha_{\text{drift}} = 0.25$, $\gamma_{\text{drift}} = 2$, \\
      cut-off function
      $
      \varphi_{\text{drift}}(\ell)=
    \begin{cases}
     0, & \ell \le 1, \\
     \ell, & \ell \ge 2, \\
     \ell^3 - 2\ell^2 + \ell, & \ell \in (1,2).
     \end{cases}
      $\\
      Exponential regularization strategy in Eq.~\eqref{3-4}: \\
      $\alpha_{\text{diff}} = 0.2$, $\gamma_{\text{diff}} = 1.5$, $\beta_{\text{diff}} = 0.2$.\\
    \end{tabular} \\
    \bottomrule
  \end{tabular}
  \label{para}
\end{table}
\par The initial velocity is generated from a bimodal density
\begin{align*}
f_0(v) &= \beta_1 \left[ e^{-\frac{(v-\mu)^2}{2\sigma^2}} + e^{-\frac{(v+\mu)^2}{2\sigma^2}} \right]
\end{align*}
with $\sigma^2 = 0.1$, $\mu = 0.25$, and $\beta_1$ such that $\int_{\mathbb{R}} f_0(v) dv = 1$. The initial position is set as $X(0)=0$.
\par The error metrics used in the  numerical tests are now introduced, taking $V$-errors as examples. For the $X$-errors, it suffices to replace all occurrences of $V$ with $X$.
\par \textbf{Notations}: randomly select $N_z=30$ values of $z$, denoted as $\mathcal{Z} = \{z_1, z_2, \dots, z_{N_z}\}$, and randomly choose $N_v=10$ initial velocities, i.e., $N_v$ particles out of the total $N=1000$ particles, denoted as $\mathcal{V} = \{v_{0,1}, v_{0,2}, \dots, v_{0,N_v}\}$. The set $\mathcal{V}$ is determined in advance and kept unchanged for all $z \in \mathcal{Z}$. Let $V_{\text{exact}}(z, v_{0,j},t)$ and $V_{\text{num}}(z, v_{0,j},t)$ denote the exact and numerical solutions corresponding to $z$ and the initial velocity $v_{0,j}$ respectively.

\par \textbf{Error for a single $z$ value}: For any $z \in \mathcal{Z}$, the error $E(z)$ is defined as follows. For each initial velocity $v_{0,j} \in \mathcal{V}$, compute the absolute difference between the exact  and  numerical solution at the final time $t=1$, and then average over the $N_v$ particles in $\mathcal{V}$, i.e.,
\begin{align*}
E(z) = \frac{1}{N_v} \sum_{j=1}^{N_v}  \left( \left| V_{\text{exact}}(z, v_{0,j},1) - V_{\text{num}}(z, v_{0,j},1) \right| \right).
\end{align*}
\par The $V$-error is then defined as the average of $E(z)$ for the $N_z$ values of $z$ in $\mathcal{Z} $, i.e.,
\begin{align}
\mathcal{E}(V) = \frac{1}{N_z} \sum_{i=1}^{N_z} E(z_i).\label{err1}
\end{align}
\par We also define 
\begin{align*}
 E_{\max}(z) = \frac{1}{N_v} \sum_{j=1}^{N_v} \left( \max_{t \in [0,1]} \left| V_{\text{exact}}(z, v_{0,j},t) - V_{\text{num}}(z, v_{0,j},t) \right| \right),
\end{align*}
and
\begin{align}
 \mathcal{E}_{\max}(V) = \frac{1}{N_z} \sum_{i=1}^{N_z} E_{\max}(z_i). \label{err2}  
\end{align}
The error metrics \eqref{err1}-\eqref{err2} are used  in Sections 5.4-5.6. In Section 5.3, a variant of these metrics will be employed to verify the weak convergence of the str-RB algorithm.
\par To facilitate the derivation of exact solutions for comparison, $P(z,\cdot)=P(z)$ is assumed in most of the following tests, i.e., $P$ depends only on the random variable $z$. For a general $P(z,\cdot)$, however, the str-RB algorithm follows the same implementation procedure, since the main difficulty in simulating of these particle systems lies in the singularities of $P$ and $D$ with respect to the random variable $z$. Therefore, simplifying $P(z,\cdot)$ to $P(z)$ in the majority of the tests does not affect the credibility of the assessment of the algorithm. To further support this claim, a test where $P$ depends on both $z$ and $X$ is provided in Test 4 of Section 5.4.
\par With $P(z,\cdot)$ simplified to $P(z)$, the exact solutions are given by
\begin{align*}
V_i(t) &= V_i(0) e^{-P(z) t} + \frac{1}{N} \sum_{k=1}^N V_k(0) \left( 1 - e^{-P(z) t} \right) \\
&+ \sqrt{2D(z)} \int\nolimits_{0}^{t} e^{-P(z) (t-y)} dW_i(y) + \sqrt{2D(z)} \int\nolimits_{0}^{t} \left( 1 - e^{-P(z) (t-y)} \right) d\bar{W}(y),
\end{align*}
and
\begin{align*}
X_i(t) &= X_i(0) + \frac{V_i(0)}{P(z)} \left( 1 - e^{-P(z) t} \right) \\
&+ \frac{t}{N} \sum_{k=1}^N V_k(0) - \frac{1}{N P(z)} \sum_{k=1}^N V_k(0) + \frac{1}{N P(z)} \sum_{k=1}^N V_k(0) e^{-P(z) t} \\
&+ \frac{\sqrt{2D(z)}}{P(z)} \int\nolimits_{0}^{t} \left( 1 - e^{-P(z) (t-u)} \right) dW_i(u) \\
&+ \sqrt{2D(z)} \int\nolimits_{0}^{t} (t-u) d\bar{W}(u) - \frac{\sqrt{2D(z)}}{P(z)} \int\nolimits_{0}^{t} \left( 1 - e^{-P(z) (t-u)} \right) d\bar{W}(u),
\end{align*}
where $\displaystyle\bar{W}(y)=\displaystyle\frac{1}{N} \sum_{i=1}^{N} W_i(y)$.
It can be  observed from the exact solutions that the singularities in $D(z)$ are the primary source of numerical instability. The remainder of this section therefore concentrates on tests where $D(z)$ has singularities.
\subsection{Expectation and variance of the regularized diffusion term}
\par To validate the asymptotic behavior established in Proposition \ref{proposition:model_perturbation}, we examine numerically how the expectation and variance of the regularized diffusion term $\tilde{D}$ vary with the time step $\tau$. The tests are conducted with time steps $\tau=2^{-5}, 2^{-6}, \dots, 2^{-10}$, and the parameters $\alpha_{\text{diff}} = 0.2,~\beta_{\text{diff}} = 0.2,~\text{and}~\gamma_{\text{diff}}=1$.
\par
\begin{figure}[htbp]
    \centering
    \includegraphics[width=0.48\textwidth]{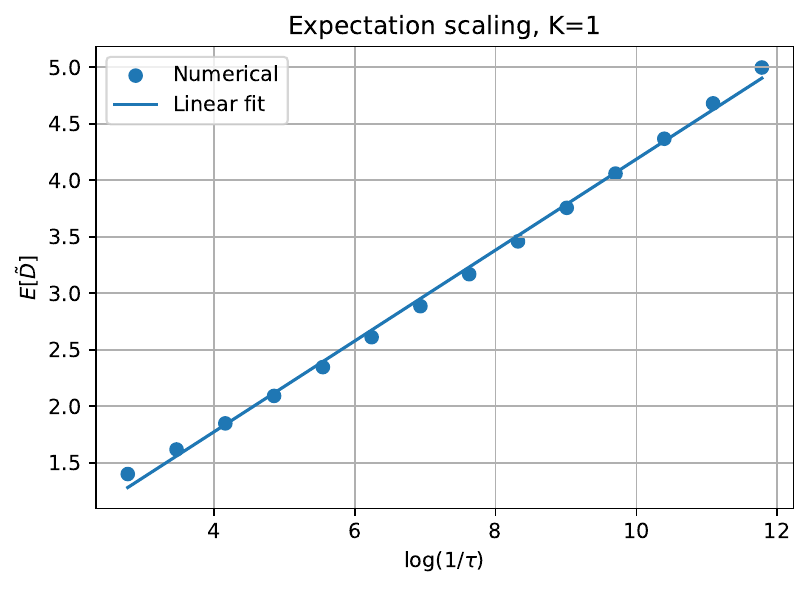}
    \includegraphics[width=0.48\textwidth]{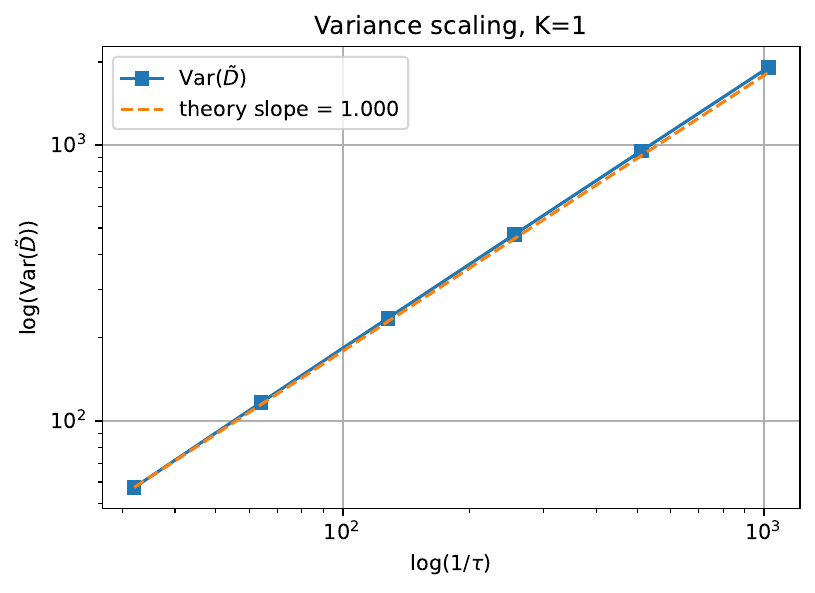}
    \caption{Variation of $\mathbb{E}[\tilde{D}]$ and $\mathrm{Var}(\tilde{D})$ with respect to $\tau$ for $K=1$.}
    \label{K=1}
\end{figure}

\begin{figure}[htbp]
    \centering
    \includegraphics[width=0.48\textwidth]{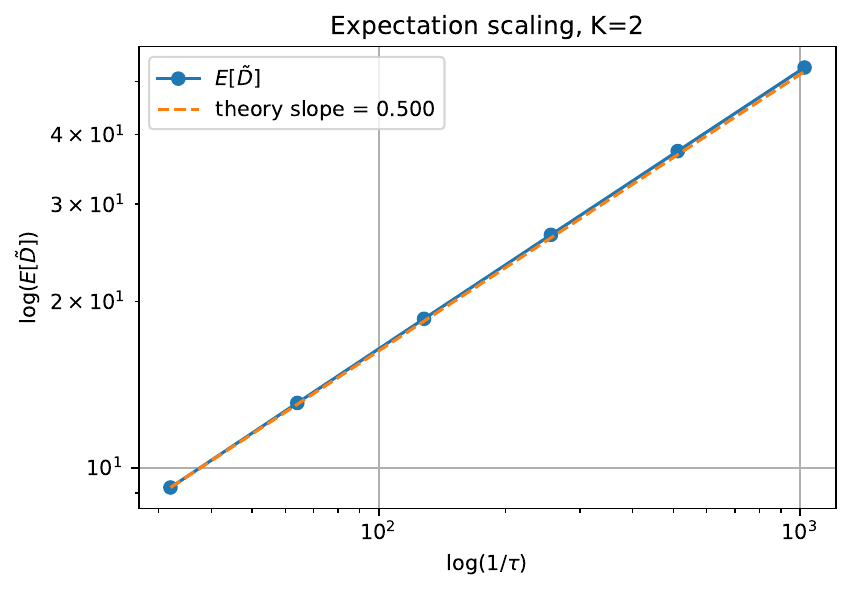}
    \includegraphics[width=0.48\textwidth]{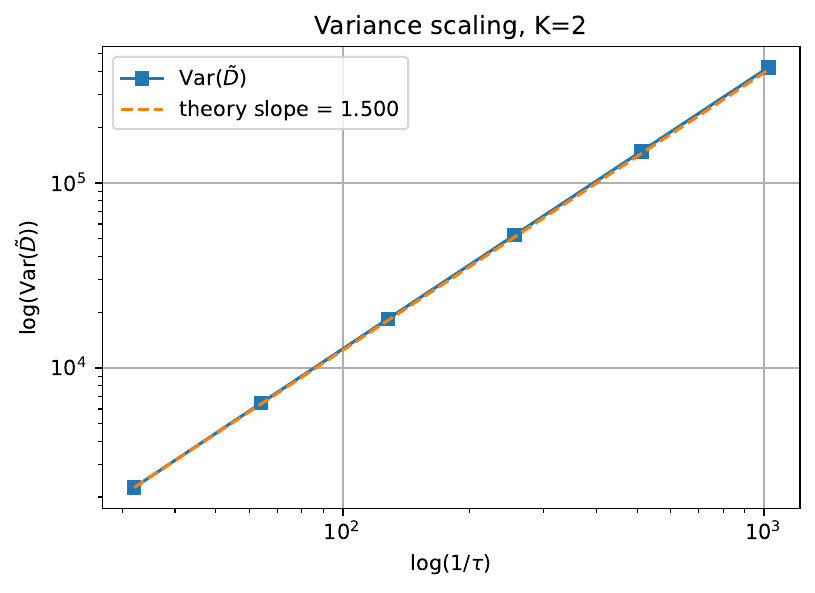}
    \caption{Variation of $\mathbb{E}[\tilde{D}]$ and $\mathrm{Var}(\tilde{D})$ with respect to $\tau$ for $K=2$.}
    \label{K=2}
\end{figure}

 Figure \ref{K=1} presents the results for the singularity order $K=1$, i.e., $D(z)=\frac{1}{1-z}$, while Figure \ref{K=2} shows the corresponding results for $K=2$. For the expectation in the $K=1$ case, the horizontal axis is $\log(1/\tau)$ and the vertical axis is $\mathbb{E}[\tilde{D}]$, so that the predicted logarithmic growth appears as a straight line; the other three plots are presented on log-log scales to demonstrate their power-law behavior.

\par The numerical results are consistent with the asymptotic behavior characterized in Proposition \ref{proposition:model_perturbation}. In particular, for $K=1$, $\mathbb{E}[\tilde{D}]$ exhibits an approximately linear dependence on $\log(1/\tau)$, while the slopes observed in the remaining log-log plots are consistent with the theoretical rates derived in Eqs. \eqref{17} and \eqref{18}.
\subsection{Convergence test}
\par Having verified the asymptotic behavior of the regularized diffusion coefficient, we now turn to the convergence of the str-RB algorithm. In view of the modeling perturbation analyzed in Section 4.3 and the numerical error discussed in Section 4.4, the error measured in the experiments contains both the numerical error of the regularized system and the perturbation introduced by the regularization strategy. Consequently, the observed convergence behavior in the numerical tests reflects the combined effect of both error sources, rather than the effect  of the regularized system alone. Since the regularization-induced perturbation is characterized statistically, weak errors are adopted as the convergence metrics in the following tests.
\par We first introduce the weak error metrics. Using the notations in Section 5.1, and again taking the $V$ errors as examples, we define for each $z$
\begin{align*}
\tilde{E}(z) = \frac{1}{N_v} \sum_{j=1}^{N_v}  \left( V_{\text{exact}}(z, v_{0,j},1) - V_{\text{num}}(z, v_{0,j},1)  \right),
\end{align*}
with $N_v=10$, the same as in Section 5.1. Similar to Section 5.1, the error of $V$ is defined as the average of $\tilde{E}(z)$ for the $N_z=30$ values of $z$ in $\mathcal{Z} $, i.e.,
\begin{align*}
\tilde{\mathcal{E}}(V) = \frac{1}{N_z} \sum_{i=1}^{N_z} \tilde{E}(z_i).
\end{align*}
\par The convergence tests are performed with $N_s=200$ random seeds for each time step, where the $N_s$ seeds are fixed in advance and independent of the step size. For each seed, the corresponding errors $\tilde{\mathcal{E}}_{k}(V)$ are computed, and  the final weak error $\mathcal{{E}}_{\text{weak}}(V)$ is defined as the absolute value of the average over the seeds, i.e.,
\begin{align}
  \mathcal{{E}}_{\text{weak}}(V)=\left|\frac{1}{N_s}\sum_{k=1}^{N_s}  \tilde{\mathcal{E}}_{k}(V)\right|,\label{err3}
\end{align}
which serves as an estimator for the absolute value of the expectation of the error. 
\par In the convergence tests, we run the algorithm with time steps $\tau = 1/64, 1/128, 1/256$, and $1/512$ on a highly singular test
\begin{align*}
P(z)=D(z)=\frac{1}{(1-z)^4},
\end{align*}
with $\gamma_{\text{drift}}=1.5,~\gamma_{\text{diff}}=3.5,~\text{and}~\beta_{\text{diff}}=1/4$ for str-RB.
\par
\begin{figure}[htbp]
  \centering
  \includegraphics[width=0.48\textwidth]{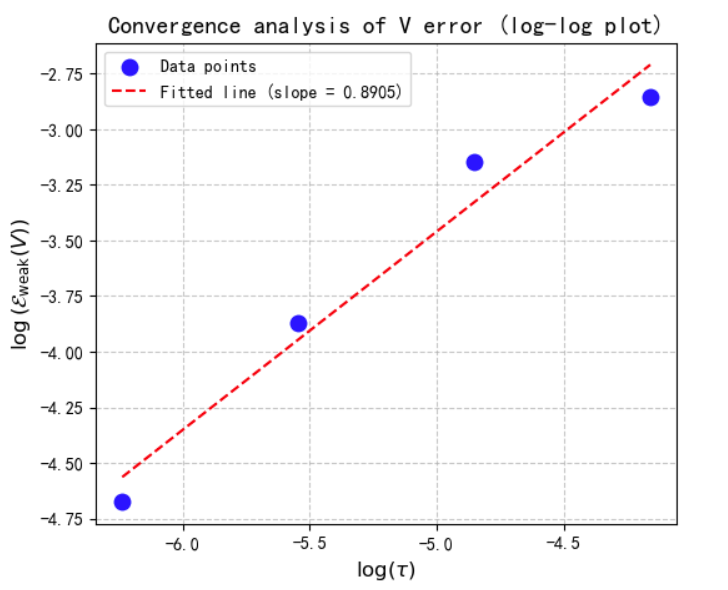} 
  \includegraphics[width=0.48\textwidth]{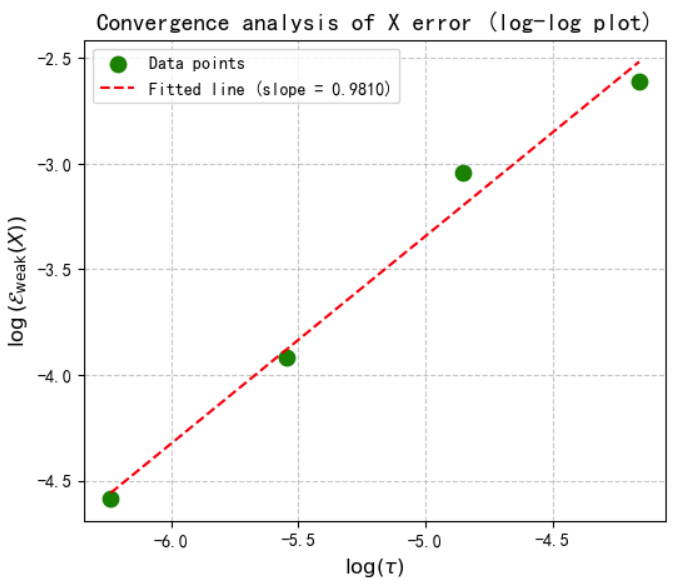} 
  \caption{Weak convergence of V-error and X-error}
  \label{fig:weak_convergence1}
\end{figure}

 The numerical results for different time steps are summarized in Table \ref{tab:weak_convergence}.  Figure \ref{fig:weak_convergence1} displays the weak errors against the time step in a log-log scale. The fitted slopes are approximately 0.89 for $\mathcal{{E}}_{\text{weak}}(V)$ and 0.98 for $\mathcal{{E}}_{\text{weak}}(X)$, indicating a convergent behavior of the str-RB algorithm. A rigorous theoretical proof of the convergence order for the str-RB algorithm in this singular setting involves substantial technical difficulties and is left for future work.
\begin{table}[htbp]
  \centering
  \caption{Weak errors at different time steps}
  \label{tab:weak_convergence}
  \begin{tabular}{lcccc}
    \toprule
    & $\tau=1/64$ & $\tau=1/128$ & $\tau=1/256$ & $\tau=1/512$ \\
    \midrule
    $\mathcal{E}_{\text{weak}}(V)$ & 0.05745 & 0.04290 & 0.02089 & 0.00933 \\
    $\mathcal{E}_{\text{weak}}(X)$ & 0.07349 & 0.04773 & 0.01995 & 0.01019 \\
    \bottomrule
  \end{tabular}
\end{table}
\subsection{Numerical tests with high-order singularities} 
\par This section presents detailed numerical tests for highly singular cases, where the diffusion term $D(z)$ exhibits singularities of order 3-4. These cases are not only the main focus of this work, but also those in which the str-RB algorithm demonstrates its key advantages over existing methods in enhancing stability and reducing errors.
\par \textbf{Test 1}:
We start with a test where $P(z)$ exhibits strong fluctuations, and $D(z)$ has singularities of order 3. Take 
\begin{align*}
P(z)=1+\sin(100\pi z),~D(z)=\frac{1}{(1-z)^3},    
\end{align*}
with $\beta_{\text{diff}}=1/4$ for str-RB. Since $D(z)$ exhibits strong singularities and $P(z)$ attains small values at many points, $\dfrac{\sqrt{D(z)}}{P(z)}$ can become very large. Accordingly, a larger $\beta_{\text{diff}}$ is chosen to strengthen the regularization effect.
\par
\begin{figure}[htbp]
  \centering
  \includegraphics[width=0.8\textwidth]{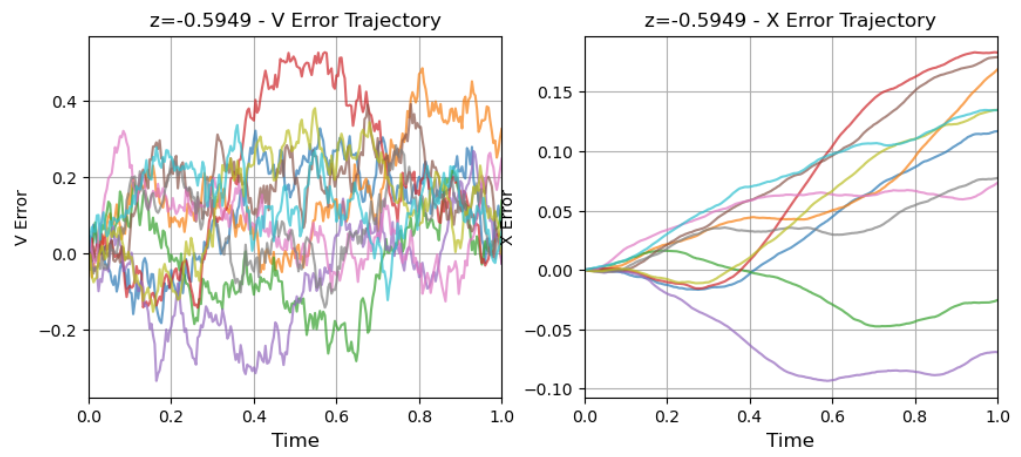} 
  \caption{$z=-0.5949$ error trajectories for Test 1}
  \label{figure1}
\end{figure}

 The results of this test are presented in Table \ref{table1}, with the errors defined in Eqs. \eqref{err1} and \eqref{err2}. The str-RB error trajectories of ten particles for $z=-0.5949$ are illustrated in Figure \ref{figure1}. The vertical axes, labeled $V$-Error and $X$-Error, represent $V_{\text{exact}}-V_{\text{num}}$ and $X_{\text{exact}}-X_{\text{num}}$ respectively. In this test, the high-order singularity in the diffusion term causes existing s-RB and st-RB algorithms to produce errors exceeding 300, while the errors of the str-RB algorithm remain below 5.
\begin{table}[htbp]
  \centering
  \caption{Results of Test 1}
  \label{table1}
  \begin{tabular}{lcccc}
    \toprule
    &s-RB & st-RB & stt-RB & str-RB \\
    \midrule
    $\mathcal{E}(V)$ & 967.12 & 1488.87 & 1.5539 & 1.3880 \\
    $\mathcal{E}(X)$ & 585.29 & 796.89 & 0.9421 & 0.8371 \\
    $\mathcal{E}_{\max}(V)$ & 2476.46 & 2799.67 & 4.1984 & 3.7292 \\
    $\mathcal{E}_{\max}(X)$ & 651.34 & 845.04 & 1.0666 & 0.9472\\
    \bottomrule
  \end{tabular}
\end{table}
\par \textbf{Test 2}: In this test, the order of singularity in $D(z)$ is increased from 3 to 4. Consider
\begin{align*}
P(z) = 
\begin{cases} 
1 + \sin z, & z \in [-1, -0.5), \\
z^2, & z \in [-0.5, 0.5), \\
\exp(z), & z \in [0.5, 1],\\
\end{cases}
~~~D(z)=\frac{1}{(1-z)^4},
\end{align*}
 where $P(z)$ has discontinuities, and $D(z)$ has a singularity of order 4. In the taming function, $\gamma_{\text{drift}}$ is adjusted to 1 for str-RB, since the non-smoothness of $P(z)$ is relatively mild and a weaker taming strength suffices.
\par
\begin{figure}[htbp]
  \centering
  \includegraphics[width=0.8\textwidth]{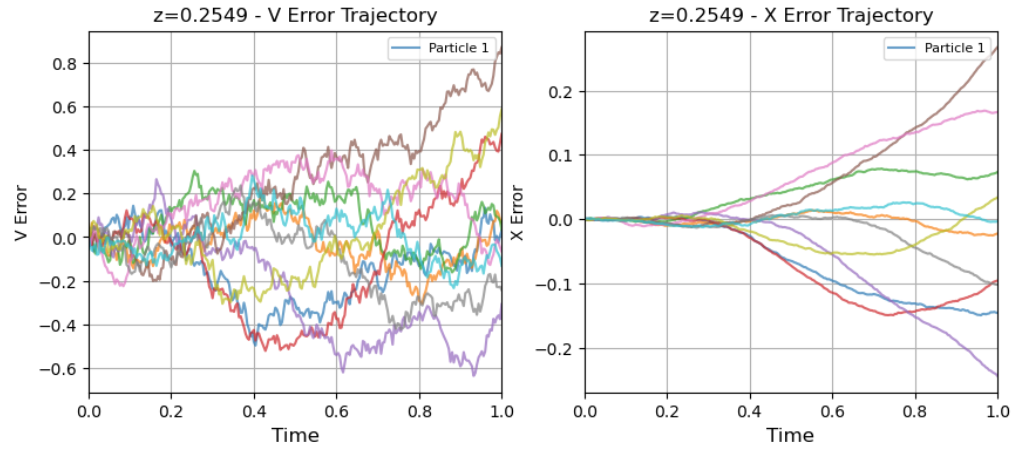} 
  \caption{$z=0.2549$ error trajectories for Test 2}
  \label{figure2}
\end{figure}

 Table \ref{table2} summarizes the results of this test, while Figure \ref{figure2} depicts the str-RB error trajectories of ten particles for $z=0.2549$. The diffusion term now has a fourth‑order singularity,  higher than that in Test 1, so the existing methods produce errors over $10^4$, but the errors of the str-RB algorithm are still controlled within 15.
\begin{table}[htbp]
  \centering
  \caption{Results of Test 2}
  \label{table2}
  \begin{tabular}{lcccc}
    \toprule
    &s-RB & st-RB & stt-RB & str-RB \\
    \midrule
    $\mathcal{E}(V)$ & 184081 & 310122 & 4.63 & 4.28 \\
    $\mathcal{E}(X)$ & 101079 & 166007 & 2.47 & 2.18 \\
    $\mathcal{E}_{\max}(V)$ & 484397 & 583291 & 13.27 & 12.16\\
    $\mathcal{E}_{\max}(X)$ & 116131 & 176044 & 2.93 & 2.62 \\
    \bottomrule
  \end{tabular}
\end{table}
\par \textbf{Test 3}: We further consider a case where both $P(z)$ and $D(z)$ have fourth-order singularities:
\begin{align*}
P(z)=D(z)=\frac{1}{(1-z)^4}.
\end{align*}
 For str-RB, $\gamma_{\text{drift}}=1.5,~\gamma_{\text{diff}}=3.5,~\beta_{\text{diff}}=1/4$ is chosen since this test is highly singular and a stronger regularization is required.
\par
\begin{figure}[htbp]
  \centering
  \includegraphics[width=0.8\textwidth]{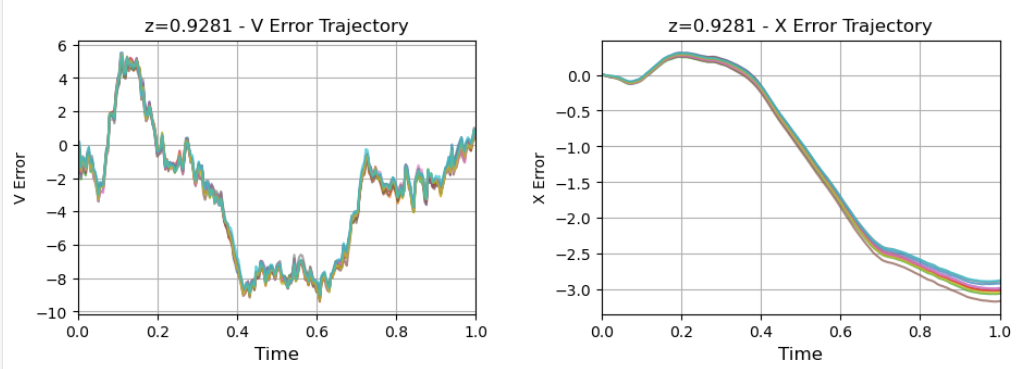} 
  \caption{$z=0.9281$ error trajectories for Test 3}
  \label{figure3}
\end{figure}

 The numerical results of this test are provided in Table \ref{table3}, with the str-RB error trajectories of ten particles at $z=0.9281$ shown in Figure \ref{figure3}. With fourth‑order singularities in both the drift and diffusion terms, the existing methods suffer from more serious error blow-up, while the str-RB algorithm still ensures stability and accuracy.
\begin{table}[htbp]
  \centering
  \caption{Results of Test 3}
  \label{table3}
  \begin{tabular}{lcccc}
    \toprule
    &s-RB & st-RB & stt-RB & str-RB \\
    \midrule
    $\mathcal{E}(V)$ & n & 310125 & 0.5329 & 0.4441 \\
    $\mathcal{E}(X)$ & n & 166010 & 0.3889 & 0.3260 \\
    $\mathcal{E}_{\max}(V)$ & n & 583304 & 1.4205 & 1.4039 \\
    $\mathcal{E}_{\max}(X)$ & n & 176047 & 0.4035 & 0.3495 \\
    \bottomrule
  \end{tabular}
  \\
  \centering      
  \parbox{\linewidth}{\centering\footnotesize (``n'' refers to errors exceeding $10^{100}$)}
\end{table}
\par \textbf{Test 4}: To examine the algorithm in a more general setting, $P$ is now allowed to depend on both $z$ and $X$. We take
\begin{align*}
P(z,X_i,X_j)=\frac{1}{(1-z)^4}\cdot\exp(|X_i-X_j|),~D(z)=\frac{1}{(1-z)^4},
\end{align*}
 where both $P(z)$ and $D(z)$ have singularities of order 4 and $P$ depends on $z$, $X_i$, and $X_j$. In this test, the singularities of $P$ and $D$ are of the same order as those in Test 3, so identical parameters are used for str-RB: $\gamma_{\text{drift}}=1.5,~\gamma_{\text{diff}}=3.5,~\beta_{\text{diff}}=1/4$. 
\par
\begin{figure}[htbp]
  \centering
  \includegraphics[width=0.8\textwidth]{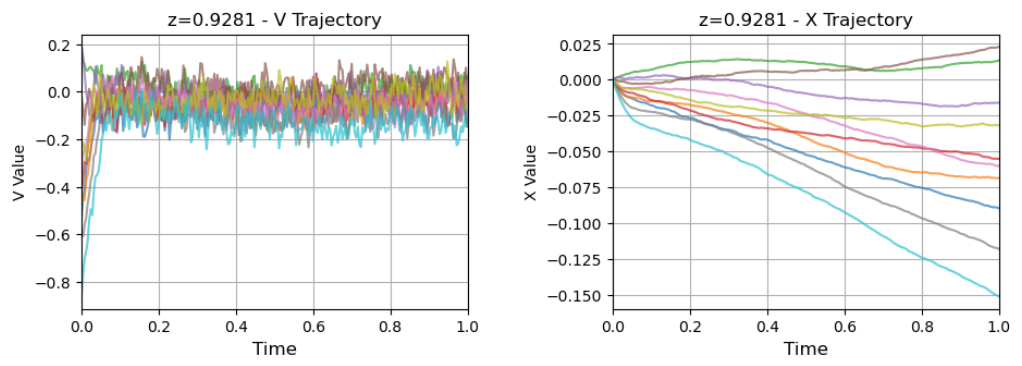} 
  \caption{$z=0.9281$ trajectories for Test 4}
  \label{fig:test4}
\end{figure}

 The trajectories of ten particles at $z=0.9281$ are shown in Figure \ref{fig:test4}. In this case, as the exact solution is difficult to obtain analytically, the trajectories of $V$ and $X$, rather than the errors, are presented. This test shows that when the coefficients are singular and $P$ depends on both $z$ and $X$, our str-RB algorithm still maintains numerical stability without blow-up.
\par In Tests 1-3 above, the str-RB algorithm reduces errors by at least two orders of magnitude compared to the s-RB and st-RB algorithms, and outperforms the enhanced baseline stt-RB algorithm. Furthermore, the str-RB algorithm does not encounter numerical blow-up throughout the evolution of the velocity and position. In all three tests, $D(z)$ has a singularity of order at least 3, while $P(z)$ is discontinuous, highly oscillatory, or singular. Under these challenging conditions, the str-RB algorithm consistently outperforms both the existing methods and the enhanced baseline, illustrating its robustness and accuracy. Test 4 further demonstrates that, when the coefficients are singular and $P$ depends on both $z$ and $X$, our algorithm still maintains good numerical stability.
\subsection{Additional numerical tests}
\par We then provide some additional numerical tests to further verify the effectiveness of the str-RB algorithm, and the results are presented in Tables \ref{table4}-\ref{table5}. In these tests, the choice of parameters is primarily based on the strength of the singularity. In general, stronger singularities require stronger taming and regularization.
\begin{table}[htbp]
  \centering
  \caption{$\mathcal{E}(V)$ of additional tests}
  \small
  \label{table4}
  \begin{tabular}{lcccc}
    \toprule
    Test case &s-RB & st-RB & stt-RB & str-RB \\
    \midrule
    $P(z)=\sin(z)+2,D(z)=\exp(z)$, & & & & \\
    st-RB: $\gamma_{\text{drift}}=0.1$; stt-RB: $\gamma_{\text{drift}}=\gamma_{\text{diff}}=0.1$; & & & & \\
    str-RB: $\gamma_{\text{drift}}=0$, $\gamma_{\text{diff}}=1$ & 0.0199 & 0.0199 & 0.0199 & 0.0204 \\
    \midrule
    $P(z)=1+\sin(100\pi z)$, $D(z)=1/(1-z)^2$, & & & & \\
    str-RB: $\beta_{\text{diff}}=1/4$ & 6.7787 & 10.1327 & 0.4767 & 0.4263 \\
    \midrule
    $P(z)=\exp(2z)/2$, $D(z)=1/(1-z)^2$, & & & & \\
    str-RB: $\gamma_{\text{drift}}=2.5,~\beta_{\text{diff}}=1/4$ & 5.6352 & 10.2803 & 0.3554 & 0.3429 \\
    \midrule
    $P(z)=
    \begin{cases}
    1+\sin z, & z\in[-1,-0.5), \\
    z^2, & z\in[-0.5,0.5), \\
    \exp(z), & z\in[0.5,1];
    \end{cases}
    ~D(z)=1/(1-z)^2$, & & & & \\
    str-RB: $\gamma_{\text{drift}}=3.5$, $ \beta_{\text{diff}}=1/4$ & 6.1078 & 10.2411 & 0.4427 & 0.4298 \\
    \midrule
    $P(z)=
    \begin{cases}
    1+\sin z, & z\in[-1,-0.5), \\
    z^2, & z\in[-0.5,0.5), \\
    \exp(z), & z\in[0.5,1];
    \end{cases}
    ~D(z)=1/(1-z)^3$ & & & & \\
    str-RB: $\beta_{\text{diff}}=1/4$ & 883.36 & 1489.31 & 1.4158 & 1.2702\\
    \midrule
    $P(z)=1$, $D(z)=1/(1-z)^4$, & & & & \\
    str-RB: $\gamma_{\text{drift}}=0.5$ & 198668 & 310119 & 5.71 & 5.29 \\
    \midrule
    $P(z)=1+\sin(100\pi z)$, $D(z)=1/(1-z)^4$, & & & & \\
    str-RB: $\gamma_{\text{diff}}=2$ & 201465 & 310120 & 5.19 & 4.87 \\
    \midrule
    $P(z)=D(z)=1/(1-z)$, & & & & \\
    str-RB: $\gamma_{\text{drift}}=0,\gamma_{\text{diff}}=0.5$ & 0.1329 & 0.2946 & 0.2221 & 0.1346 \\
    \midrule
    $P(z)=1/(1-z)$, $D(z)=1/(1-z)^2$, & & & & \\
    str-RB: $\beta_{\text{diff}}=1/4$ & 2.1687 & 10.3797 & 0.4388 & 0.3609 \\
    \midrule
    $P(z)=D(z)=1/(1-z)^2$, & & & & \\
    str-RB: $\gamma_{\text{diff}}=2.5,\beta_{\text{diff}}=1/4$ & n & 10.4999 & 0.3826 & 0.3131 \\
    \midrule
    $P(z)=1/(1-z)^3$, $D(z)=1/(1-z)^2$, & & & & \\
    str-RB: $\gamma_{\text{diff}}=3,\beta_{\text{diff}}=1/4$ & n & 10.5436 & 0.4092 & 0.3148 \\
    \midrule
    $P(z)=1/(1-z)^2$, $D(z)=1/(1-z)^3$, & & & & \\
    str-RB:  $\gamma_{\text{drift}}=\gamma_{\text{diff}}=2.5,\beta_{\text{diff}}=1/4$ & n & 1490.29 & 0.6166 & 0.5509 \\
    \midrule
    $P(z)=D(z)=1/(1-z)^3$, & & & & \\
    str-RB: $\gamma_{\text{diff}}=3.5,\beta_{\text{diff}}=1/4$ & n & 1490.37 & 0.4608 & 0.3717 \\
    \bottomrule
  \end{tabular}
  \\ 
  \centering      
  \parbox{\linewidth}{\centering\footnotesize (``n'' refers to errors exceeding $10^{100}$)}
\end{table}
\begin{table}[htbp]
  \centering
  \caption{$\mathcal{E}(X)$ of additional tests}
  \small
  \label{table5}
  \begin{tabular}{lcccc}
    \toprule
    Test case &s-RB & st-RB & stt-RB & str-RB \\
    \midrule
    $P(z)=\sin(z)+2,D(z)=\exp(z)$, & & & & \\
    st-RB: $\gamma_{\text{drift}}=0.1$; stt-RB: $\gamma_{\text{drift}}=\gamma_{\text{diff}}=0.1$; & & & & \\
    str-RB: $\gamma_{\text{drift}}=0$, $\gamma_{\text{diff}}=1$ & 0.0094 & 0.0094 & 0.0094 & 0.0091 \\
    \midrule
    $P(z)=1+\sin(100\pi z)$, $D(z)=1/(1-z)^2$, & & & & \\
    str-RB: $\beta_{\text{diff}}=1/4$ & 4.0870 & 5.4403 & 0.2790 & 0.2488 \\
    \midrule
    $P(z)=\exp(2z)/2$, $D(z)=1/(1-z)^2$, & & & & \\
    str-RB: $\gamma_{\text{drift}}=2.5,~\beta_{\text{diff}}=1/4$ & 3.0442 & 5.5243 & 0.1839 & 0.1484 \\
    \midrule
    $P(z)=
    \begin{cases}
    1+\sin z, & z\in[-1,-0.5), \\
    z^2, & z\in[-0.5,0.5), \\
    \exp(z), & z\in[0.5,1]
    \end{cases},
    D(z)=1/(1-z)^2$, & & & & \\
    str-RB: $\gamma_{\text{drift}}=3.5$, $\beta_{\text{diff}}=1/4$ & 3.3505 & 5.4969 & 0.2324 & 0.2018 \\
    \midrule
    $P(z)=
    \begin{cases}
    1+\sin z, & z\in[-1,-0.5), \\
    z^2, & z\in[-0.5,0.5), \\
    \exp(z), & z\in[0.5,1]
    \end{cases},
    D(z)=1/(1-z)^3$ & & & & \\
    str-RB: $\beta_{\text{diff}}=1/4$ & 485.09 & 797.12 & 0.7567 & 0.6428 \\
     \midrule
    $P(z)=1$, $D(z)=1/(1-z)^4$, & & & & \\
    str-RB: $\gamma_{\text{drift}}=0.5$ & 120622 & 166006 & 3.52 & 3.26 \\
    \midrule
    $P(z)=1+\sin(100\pi z)$, $D(z)=1/(1-z)^4$, & & & & \\
    str-RB: $\gamma_{\text{diff}}=2$ & 121936 & 166007 & 3.18 & 2.97 \\
    \midrule
    $P(z)=D(z)=1/(1-z)$, & & & & \\
    str-RB: $\gamma_{\text{drift}}=0,\gamma_{\text{diff}}=0.5$ & 0.0620 & 0.1500 & 0.0846 & 0.0633\\
    \midrule
    $P(z)=1/(1-z)$, $D(z)=1/(1-z)^2$, & & & & \\
    str-RB: $\beta_{\text{diff}}=1/4$ & 0.2501 & 5.6243 & 0.1393 & 0.1269 \\
    \midrule
    $P(z)=D(z)=1/(1-z)^2$, & & & & \\
    str-RB: $\gamma_{\text{diff}}=2.5,\beta_{\text{diff}}=1/4$& n & 5.7111 & 0.1600 & 0.1539 \\
    \midrule
    $P(z)=1/(1-z)^3$, $D(z)=1/(1-z)^2$, & & & & \\
    str-RB: $\gamma_{\text{diff}}=3,\beta_{\text{diff}}=1/4$ & n & 5.7466 & 0.1938 & 0.1670 \\
    \midrule
    $P(z)=1/(1-z)^2$, $D(z)=1/(1-z)^3$, & & & & \\
    str-RB: $\gamma_{\text{drift}}=\gamma_{\text{diff}}=2.5,\beta_{\text{diff}}=1/4$ & n & 797.83 & 0.1892 & 0.1884 \\
    \midrule
    $P(z)=D(z)=1/(1-z)^3$, & & & & \\
    str-RB: $\gamma_{\text{diff}}=3.5,\beta_{\text{diff}}=1/4$ & n & 797.89 & 0.2285 & 0.1903 \\
    \bottomrule
  \end{tabular}
  \\
  \centering
  \parbox{\linewidth}{\centering\footnotesize (``n'' refers to errors exceeding $10^{100}$)}
\end{table}
\par These additional numerical tests cover a diverse range of cases, including a smooth case ($P(z)=\sin(z)+2, D(z)=\exp(z)$), singular cases with higher-order singularity for $P(z)$ (e.g., $P(z)=1/(1-z)^3$, $D(z)=1/(1-z)^2$), higher-order singularity for $D(z)$ (e.g., $P(z)=1/(1-z)^2$, $D(z)=1/(1-z)^3$), equal singularity orders for both (e.g., $P(z)=D(z)=1/(1-z)^3$), and mixed cases where one exhibits singularities and the other is strongly oscillatory (e.g., $P(z)=1+\sin(100\pi z)$, $D(z)=1/(1-z)^4$). In smooth and mildly singular cases, the str-RB algorithm performs comparably to the comparison methods, while in cases where the order of singularities exceeds two, the str-RB algorithm reduces the errors by factors ranging from 10 to over $10^{4}$ compared to the existing methods and performs better than the enhanced baseline, further confirming its stability and effectiveness.
\subsection{Sensitivity analysis}
\par To verify the robustness of the str-RB algorithm against variations of its core parameters, a sensitivity analysis is conducted in this section. Using Test 1 in Section 5.4, i.e., $P(z)=1+\sin(100\pi z)$ and $D(z)=\frac{1}{(1-z)^3}$,  each parameter is varied individually while all the others are kept fixed at their values in Section 5.4. The detailed results are shown in Table \ref{table6}, where for each parameter, the values adopted in Section 5.4 and the corresponding errors are also listed for comparison. A more thorough sensitivity analysis of the key parameter $\gamma_{\text{diff}}$ is conducted for values from 1.0 to 1.8 in steps of 0.05, and the associated errors $\mathcal{E}(V)$ and $\mathcal{E}(X)$ are calculated, with the results presented in Figure \ref{figure4}.
\begin{table}[htbp]
  \centering
  \caption{Sensitivity analysis results of key parameters in the str-RB algorithm}
  \label{table6}
  \begin{tabular}{lcccc}
    \toprule
    Parameter value & $\mathcal{E}(V)$ & $\mathcal{E}(X)$ & $\mathcal{E}_{\max}(V)$ & $\mathcal{E}_{\max}(X)$ \\
    \midrule
    $\gamma_{\text{diff}}=1.25$ & 1.3532 & 0.8195 & 3.6457 & 0.9255 \\
    $\gamma_{\text{diff}}=1.5$ & 1.3880 & 0.8371 & 3.7292 & 0.9472 \\
    $\gamma_{\text{diff}}=1.75$ & 1.4378 & 0.8700 & 3.8833 & 0.9849 \\
    $\beta_{\text{diff}}=1/5$ & 1.4375 & 0.8527 & 3.7945 & 0.9627 \\
    $\beta_{\text{diff}}=1/4$ & 1.3880 & 0.8371 & 3.7292 & 0.9472 \\
    $S=3$ & 1.3787 & 0.8324 & 3.7199 & 0.9493 \\
    $S=4$ & 1.3863 & 0.8397 & 3.7328 & 0.9503 \\
    $S=5$ & 1.3880 & 0.8371 & 3.7292 & 0.9472 \\
    $M=1$ & 1.5564 & 0.9448 & 4.2078 & 1.0688 \\
    $M=2$ & 1.3880 & 0.8371 & 3.7292 & 0.9472\\
    \bottomrule
  \end{tabular}
\end{table}
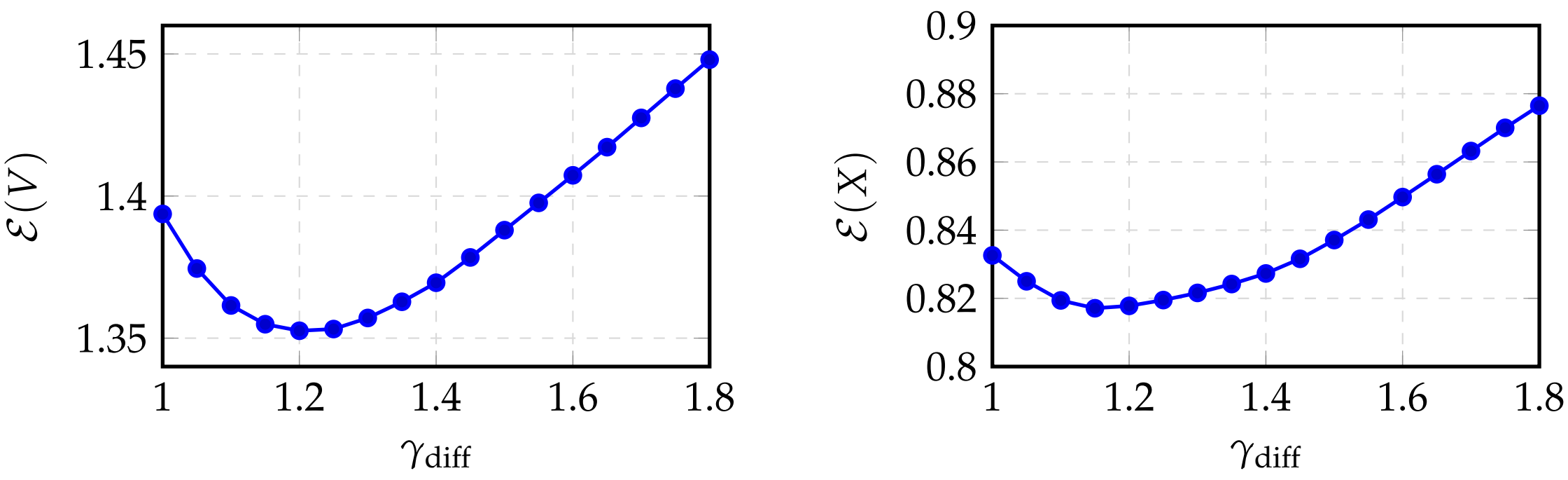
\begin{figure}[htbp]
\centering

\begin{tikzpicture}
\begin{axis}[
width=0.45\textwidth,
height=0.32\textwidth,
xlabel={$\gamma_{\text{diff}}$},
ylabel={$\mathcal{E}(V)$},
xmin=1.0, xmax=1.8,
ymin=1.34, ymax=1.46,
grid=major,
grid style={dashed,gray!30},
mark=circle*,
line width=1pt,
]
\addplot coordinates {
(1.00, 1.3937)
(1.05, 1.3745)
(1.10, 1.3615)
(1.15, 1.3549)
(1.20, 1.3526)
(1.25, 1.3532)
(1.30, 1.3571)
(1.35, 1.3628)
(1.40, 1.3695)
(1.45, 1.3784)
(1.50, 1.3880)
(1.55, 1.3976)
(1.60, 1.4073)
(1.65, 1.4172)
(1.70, 1.4275)
(1.75, 1.4378)
(1.80, 1.4480)
};
\end{axis}
\end{tikzpicture}
\hfill
\begin{tikzpicture}
\begin{axis}[
width=0.45\textwidth,
height=0.32\textwidth,
xlabel={$\gamma_{\text{diff}}$},
ylabel={$\mathcal{E}(X)$},
xmin=1.0, xmax=1.8,
ymin=0.8, ymax=0.9, 
grid=major,
grid style={dashed,gray!30},
mark=square*,
line width=1pt,
]
\addplot coordinates {
(1.00, 0.8326)
(1.05, 0.8250)
(1.10, 0.8194)
(1.15, 0.8171)
(1.20, 0.8178)
(1.25, 0.8195)
(1.30, 0.8216)
(1.35, 0.8242)
(1.40, 0.8273)
(1.45, 0.8316)
(1.50, 0.8371)
(1.55, 0.8431)
(1.60, 0.8497)
(1.65, 0.8564)
(1.70, 0.8632)
(1.75, 0.8700)
(1.80, 0.8765)
};
\end{axis}
\end{tikzpicture}

\caption{Sensitivity analysis of $\gamma_{\text{diff}}$}
\label{figure4}
\end{figure}
\par The results show that, for all parameter values in Table \ref{table6}, almost all error metrics of the str-RB algorithm are lower than the corresponding metrics of the s-RB algorithm, the st-RB algorithm and the enhanced baseline stt-RB algorithm reported in Table \ref{table1}. Figure \ref{figure4} further illustrates that, for $\gamma_{\text{diff}} \in [1, 1.8]$, the str-RB algorithm consistently outperforms the enhanced baseline: within this range, $\mathcal{E}(V)$ and $\mathcal{E}(X)$ do not exceed 1.45 and 0.88 respectively, compared to 1.5539 and 0.9421 for the stt-RB algorithm.  This further confirms that within this range, the errors are insensitive to the variation of parameters and that the str-RB algorithm maintains its superior performance over the other three algorithms. 
\subsection{Discussion of the numerical results}
\par The numerical tests in Section 5.3 demonstrate that the str-RB algorithm exhibits weak convergence. Sections 5.4 and 5.5 cover a diverse range of cases, including high-order singularities of order 3 to 4 in the diffusion term, simultaneous singularities in both drift and diffusion terms, low-order singularities of order 1 to 2, and mixed non-smoothness. The results show that the proposed str-RB algorithm outperforms the benchmark algorithms s-RB, st-RB, and stt-RB in both numerical stability and accuracy, with clear advantages in handling singularities, especially the challenging singularities in the diffusion term. Moreover, the str-RB algorithm does not exhibit numerical blow-up across all these tests. 
The existing s-RB and st-RB algorithms become unstable when $D(z)$ contains singularities of order higher than two. The enhanced baseline stt-RB, constructed in this study by extending the taming strategy to the diffusion terms, maintains stability in challenging singular problems and improves performance over s-RB and st-RB. Nevertheless, the str-RB algorithm achieves lower errors than stt-RB in all singular tests, showing better numerical stability and higher accuracy. 
\par This improved performance is attributed to the more effective treatment of singularities by the exponential regularization strategy. In cases where $D(z)$ exhibits challenging singularities, it must be regularized before performing the sG projection to avoid instability. Therefore, the s-RB and st-RB algorithms produce large errors when $D(z)$ contains singularities of order higher than two. The stt-RB algorithm tames $D(z)$ prior to the sG projection, which improves numerical stability. Nevertheless, when $D(z)$ contains high-order singularities, the algebraic damping of the taming strategy becomes insufficient, and a stronger exponential regularization is required to further enhance accuracy. Consequently, the str-RB algorithm, which applies exponential regularization to the diffusion term, achieves lower errors than the enhanced baseline stt-RB across all singular cases. Furthermore, when the key parameters vary within the considered ranges, the str-RB algorithm performs better than the comparison methods for nearly all tested values. The numerical test results confirm that the exponential regularization strategy proposed in this work is effective for a broader class of singular problems and provides a robust and accurate approach for uncertain particle systems with singular coefficients.
\section{Conclusion}
\par In this paper, we propose a novel stochastic tamed-regularized Random Batch (str-RB) algorithm for uncertain particle systems with coefficients that are singular with respect to the random variable. To address numerical blow-up caused by singularities in the diffusion term, an exponential regularization strategy is designed before the sG projection, which yields an SDE with deterministic coefficients. A random batch approach is then applied to this SDE to reduce the computational complexity, and the resulting problem is solved using a modified tamed scheme to further enhance stability. Several theoretical properties of the proposed algorithm are established, showing in particular that the exponential regularization provides stronger suppression of large diffusion coefficients while recovering the original system more accurately as the time step tends to zero. 
In a series of numerical tests, we investigate the convergence of the str-RB algorithm and compare it with the original s-RB, the st-RB and the enhanced baseline stt-RB. Numerical results show that the str-RB algorithm exhibits weak convergence, and that over a wide range of challenging singular cases, str-RB enhances stability and  yields lower errors than both existing methods and the enhanced baseline. Furthermore, the parameter sensitivity analysis verifies that str-RB maintains its accuracy and stability even when its key parameters vary within a reasonable range. The exponential regularization strategy plays a critical role in alleviating the numerical instability arising from singularities in the diffusion term together with the strong fluctuations of Brownian motion, effectively resolving a key difficulty in uncertain particle systems with challenging singular diffusion terms.
\par In the future, one may establish rigorous theoretical proofs of the convergence order of our algorithm, or develop formal selection criteria for key parameters in this algorithm. Furthermore, extending the algorithm to particle systems with similar structures as well as related multiscale problems, and combining it with other uncertainty quantification methods are also promising directions worth exploring.

\section*{Acknowledgments}
The work was partially supported by National Key Research and Development Program of China (2024YFA1016100), National Natural Science Foundation of China (12671502, 12531016) and Shanghai Jiao Tong University 2030 Initiative (WH510363004).


\end{document}